\documentclass[12pt]{article}
\usepackage{booktabs} 
\usepackage{amssymb,amsmath,fullpage,times,epsf}
\usepackage{CJK}
\usepackage{mathrsfs}
\usepackage{amssymb}
\usepackage{graphics}
\usepackage{graphicx}
\usepackage{caption}
\usepackage{graphicx}
\usepackage{float}
\usepackage{subcaption}
\usepackage{color}
\usepackage{bm}
\usepackage{amsfonts}
\usepackage{color}
\usepackage{indentfirst}
\usepackage{diagbox}
\usepackage{multirow}
\usepackage{bm}
\usepackage{cite}
\usepackage{makecell}
\usepackage{verbatim}
\usepackage{multirow,bigstrut}

\usepackage{booktabs}
\usepackage[ruled]{algorithm2e}

\begin{document}
\title{RAR-PINN algorithm for the data-driven vector-soliton solutions and parameter discovery of coupled nonlinear equations}
\author{Shu-Mei Qin$ ^{a} $, Min Li$ ^{a,} $\thanks{Corresponding author, e-mail: ml85@ncepu.edu.cn}, Tao Xu$ ^{b} $, Shao-Qun Dong$ ^{b} $\\
	\\{\em $ ^{a} $  School of Mathematics and Physics,}\\
	{\em  North China Electric Power University, Beijing 102206,
		China} \\
	{\em $ ^{b} $ College of Science, China University of Petroleum, Beijing
		102249, China} }
\date{}
\maketitle
\vspace{-5mm}
\begin{abstract}
This work aims to provide an effective deep learning framework to predict the vector-soliton solutions of the coupled nonlinear equations and their interactions. The method we propose here is a physics-informed neural network (PINN) combining with the residual-based adaptive refinement (RAR-PINN) algorithm. Different from the traditional PINN algorithm which takes points randomly, the RAR-PINN algorithm uses an adaptive point-fetching approach to improve the training efficiency for the solutions with steep gradients. A series of experiment comparisons between the RAR-PINN and traditional PINN algorithms are implemented to a coupled generalized nonlinear Schr\"{o}dinger (CGNLS) equation as an example. The results indicate that the RAR-PINN algorithm has faster convergence rate and better approximation ability, especially in modeling the shape-changing vector-soliton interactions in the coupled systems. Finally, the RAR-PINN method is applied to perform the data-driven discovery of the CGNLS equation, which shows the dispersion and nonlinear coefficients can be well approximated.
\end{abstract}


\newpage

\section{Introduction}
In recent decades, with the explosion of computing resources and available data, deep learning has been successfully applied in diverse fields, such as language understanding \cite{2011Natural}, face recognition \cite{2018Face}, risk assessment \cite{zio2018future}, mathematical physics \cite{raissi2019physics,yu2022gradient,karniadakis2021physics} and so on. Particularly, the idea of using continuous dynamical systems to model high-dimensional nonlinear functions was proposed in Ref.\cite{Weinan2017A}, i.e., the relationship between differential equations and deep learning. Furthermore, Ref.\cite{raissi2019physics} proposed the physics-informed neural network (PINN) algorithm to solve partial differential equations (PDEs). Compared with the traditional mesh-based methods, such as the finite element method and the finite difference method, PINN method is mesh-free by taking advantage of the automatic differentiation technique \cite{paszke2017automatic}. The original PINN method has been successfully applied in solving nonlinear partial differential equations with a small number of data sets \cite{haghighat2021sciann,2020Learning}, and could break the curse of dimensionality \cite{poggio2017and} in classical numerical methods. Recently, a lot of high-quality data-driven solutions for classical integrable systems were obtained by using the PINN  method \cite{sun2020bright,miao2022physics}. The breathers and rogue wave solutions of the nonlinear Schr\"{o}dinger (NLS) equation were recovered with the aid of the PINN model \cite{pu2021soliton}. What's more, a variety of femtosecond optical soliton solutions of the high-order NLS equation were investigated with the PINN method \cite{2021Data}.

In order to improve the accuracy and convergence rate of the traditional PINN (TPINN) algorithm, many modified algorithms have been proposed. For instance, Ref. \cite{wang2021understanding} introduced the gradient pathologies into PINN algorithm based on the idea of adaptive loss function weights, Ref. \cite{pu2021soliton} presented neuron-wise locally adaptive activation function to derive localized wave solutions and Ref. \cite{2020Extended} proposed a generalized space-time domain decomposition approach to solve PDEs on arbitrary complex-geometry domains. Each of these approaches makes change to the neural network itself. In this paper, we will propose a new RAR-PINN method, which combines a residual-based adaptive refinement method to improve the training efficiency. The idea is that, solutions with steep gradients can be learned better by putting more points near the sharp front. Moreover, we will compare the TPINN and RAR-PINN algorithms by comparing the residuals, loss function curves, absolute errors, and relative errors in detail.

Most of the previous works investigated the data-driven scalar solitons. In this study, we would like to explore the data-driven vector solutions of integrable coupled nonlinear equations. Specially, the coupled generalized nonlinear Schr\"{o}dinger (CGNLS) equation \cite{2010Integrable} is used as an example to illustrate, namely
\begin{equation}
\left\{
\begin{aligned}
&ih_{1t}+h_{1xx}+2( \alpha |h_{1}|^2+\beta|h_{2}|^2+ \gamma h_{1}h_{2}^{*}+ \gamma ^{*}h_{2}h_{1}^{*})h_{1}=0,x\in[-L,L],t\in[-T,T]\\\
&ih_{2t}+h_{2xx}+2(\alpha |h_{1}|^2+\beta |h_{2}|^2+\gamma h_{1}h_{2}^{*}+\gamma ^{*}h_{2}h_{1}^{*})h_{2}=0,\label{equation}
\end{aligned}
\right.
\end{equation}
where $h_{1}$ and $h_{2}$ are slowly varying pulse envelopes, $\alpha$ and
$\beta$ are real constants corresponding to the self-phase modulation and cross-phase modulation effects, respectively. Here $\gamma$ is a complex constant corresponding to four-wave mixing, the subscripts $x$ and $t$ are the partial derivatives about the scaled distance and retarded time, and $*$ represents complex conjugation. When $\alpha=\beta$ and $\gamma=0$, the above equation reduces to the Manakov system, which is a well known two-coupled NLS equation in the literature \cite{manakov1974complete}. What's more, when $\alpha=-\beta$ and $\gamma=0$, it reduces to the mixed coupled NLS equation \cite{vijayajayanthi2008bright,agalarov2015bright}. The nonlinear system (\ref{equation}) is proved to be completely integrable for arbitrary values of the system parameters $\alpha, \gamma$ and $\beta$ through the Weiss-Tabor-Carnevale algorithm \cite{lu2013painleve}. Multihump solitons were demonstrated via experiments when the self-trapped incoherent wave packets propagate in a dispersive nonlinear medium \cite{1997Self,1998Observation}. The non-degenerate vector-soliton solutions can be applied in the fields of Bose-Einstein condensates \cite{2016Solitons,stamper1998optical} and optics communication \cite{sun2020bright}. Up until now, the Manakov system has been solved by a pre-fixed multi-stage training algorithm \cite{mo2022data}. In this paper, we will apply the RAR-PINN method to model various types of vector solitons and their interactions of Eq.~(\ref{equation}) which is general than the Manakov system.

The structure of this paper is as follows. In Section 2, the algorithm of the RAR-PINN is provided for the $N$-component coupled nonlinear equations. In Section 3, the vector one-, two- and three-soliton solutions of Eq.~(\ref{equation}) are predicted by using the TPINN and RAR-PINN methods, respectively. The relative errors, loss functions, and absolute errors for those two methods are also compared. In Section 4, we discuss how to apply the RAR-PINN method to the data-driven parameters discovery for Eq.~(\ref{equation}). Finally, conclusions are given in Section 5.

\section{Method}

Multi-layer feedforward neural networks have strong expressive power. Theoretically, as long as the number of neurons is large enough, the feedforward neural networks can simultaneously approximate any function and their partial derivatives \cite{1999Approximation}. TPINN \cite{raissi2019physics} is a class of deep learning algorithm where the PDE is satisfied through the loss function of the multi-layer neural networks. The TPINN algorithm is simple, which may not be efficient for some PDEs that exhibit solutions with steep gradients \cite{lu2021deepxde}. In order to improve the training efficiency of TPINN, the RAR-PINN method is proposed. The improved algorithm has great potential to be applied to various fields where adaptivity is needed \cite{hanna2021residual}.

In this paper, we consider the following form of $(1+1)$-dimensional $N$-coupled nonlinear time-dependent system in complex space,
\begin{equation}
	\begin{aligned}
			&i\bm{H}_t+\bm{\chi}(\bm{\lambda},\bm{H},\bm{H}_{xx},\bm{H}_{xxx},...)=0,x\in\Omega,t\in[-T,T]\\
			&\mathcal{B}(\bm{H},x,t)=0,x \in \partial \Omega,
			\label{eq1}
	\end{aligned}	
\end{equation}
and solve their soliton solutions $ \bm{H}=[h_{1}(x,t),h_{2}(x,t),...,h_{n}(x,t)] $ is a complex-valued vector. $ \bm{\chi}=[\chi_{1},\chi_{2},...,\chi_{n}] $ are nonlinear differential operators, specially, they can be taxonomized into a hyperbolic \cite{iwasaki1985cauchy}, parabolic \cite{2009The} or elliptic differential operator \cite{2010Fundamentals}. $ \bm{\lambda}=[\lambda_{1},\lambda_{2},...] $ are parameters of the equation. $ \mathcal{B}(\bm{H},x,t) $ could be  hyperbolic Neumann boundary conditions \cite{mcfall2009artificial}, Dirichlet boundary conditions \cite{heydari2014legendre},  periodic boundary conditions \cite{makov1995periodic} or mixed boundary conditions \cite{cano2018linear}. Besides, the initial condition can be simply treated as a special type of Dirichlet boundary condition on the spatia-temporal domain \cite{miao2022physics}. The subscripts $t$ and $x$ are the partial derivatives and correspond to the retarded time and scaled distance. Then, Eq.~\eqref{eq1} can be converted into
\begin{equation}
    \left\{	
	\begin{aligned}
		&ih_{1t}+\chi_{1}(\lambda,h_{1},h_{2},h_{1x},h_{2x},...)=0,\\
		&ih_{2t}+\chi_{2}(\lambda,h_{1},h_{2},h_{1x},h_{2x},...)=0,\\
		&\vdots\\
		&ih_{nt}+\chi_{n}(\lambda,h_{1},h_{2},h_{1x},h_{2x},...)=0.
		\label{eq2}
	\end{aligned}
   \right.
\end{equation}
Next, we decompose $h_{i}(x,t)$ into the real part $u_{i}(x,t)$ and the imaginary part $v_{i}(x,t)(i=1,2,...,n)$, i.e., $h_{1}(x,t)=u_{1}+v_{1}i$, $h_{2}(x,t)=u_{2}+v_{2}i,...,h_{n}(x,t)=u_{n}+v_{n}i$. Then Eq.~\eqref{eq2} reduces to
\begin{equation}
	\left\{
	\begin{aligned}
		&u_{1t}+\chi_{1u}(\lambda,u_{1},v_{1},u_{2},v_{2},...)=0,\\
		&v_{1t}+\chi_{1v}(\lambda,u_{1},v_{1},u_{2},v_{2},...)=0,\\
		&u_{2t}+\chi_{2u}(\lambda,u_{1},v_{1},u_{2},v_{2},...)=0,\\
		&v_{2t}+\chi_{2v}(\lambda,u_{1},v_{1},u_{2},v_{2},...)=0,\\
		&\vdots\\
		&u_{nt}+\chi_{nu}(\lambda,u_{1},v_{1},u_{2},v_{2},...)=0,\\
		&v_{nt}+\chi_{nv}(\lambda,u_{1},v_{1},u_{2},v_{2},...)=0.
		\label{eq3}
	\end{aligned}
	\right.
\end{equation}
Specifically, $\chi_{1u},\chi_{1v},\chi_{2u},...,\chi_{nv}$ are nonlinear functions of the soliton solutions and their spatial derivatives $u_{1},v_{1},u_{2},v_{2},u_{1x},v_{1x},u_{2x},v_{2x}...$. We define the residuals $f_{1u}(x,t),f_{1v}(x,t),f_{2u}(x,t),...$ and $f_{nv}(x,t)$, respectively
\begin{equation}
	\left\{
	\begin{aligned}
		&f_{1u}:=\hat{u}_{1t}+\chi_{1u}(\lambda,\hat{u}_{1},\hat{v}_{1},\hat{u}_{2},\hat{v}_{2},...),\\
		&f_{1v}:=\hat{v}_{1t}+\chi_{1v}(\lambda,\hat{u}_{1},\hat{v}_{1},\hat{u}_{2},\hat{v}_{2},...),\\
		&f_{2u}:=\hat{u}_{2t}+\chi_{2u}(\lambda,\hat{u}_{1},\hat{v}_{1},\hat{u}_{2},\hat{v}_{2},...),\\
		&f_{2v}:=\hat{v}_{2t}+\chi_{2v}(\lambda,\hat{u}_{1},\hat{v}_{1},\hat{u}_{2},\hat{v}_{2},...),\\
		&\vdots\\
		&f_{nu}:=\hat{u}_{nt}+\chi_{nu}(\lambda,\hat{u}_{1},\hat{v}_{1},\hat{u}_{2},\hat{v}_{2},...),\\
		&f_{nv}:=\hat{v}_{nt}+\chi_{nv}(\lambda,\hat{u}_{1},\hat{v}_{1},\hat{u}_{2},\hat{v}_{2},...).
	\end{aligned}
	\right.
\end{equation}
Considering governing equations, initial condition and boundary condition, the loss function for solving Eq.~\eqref{eq2} can be expressed as
\begin{equation}
	\begin{split}
		\begin{aligned}
			Loss=loss_{0}+loss_{b}+loss_{f},\label{Loss}
		\end{aligned}
	\end{split}
\end{equation}
where $ loss_{b} $ and $ loss_{0} $ represent the mean square error of boundary and initial condition, respectively, while $ loss_{f} $ is the collocation points error. Specifically,
\begin{equation}
	\begin{split}
		\begin{aligned}
			loss_{0}=&\frac{1}{N_{0}}\sum_{i=1}^{N_{0}}(|\hat{u}_{1}(x^{i},0)-u_{10}^i|^2+|\hat{v}_{1}(x^{i},0)-v_{10}^i|^2\\
			&+|\hat{u}_{2}(x^{i},0)-u_{20}^i|^2+|\hat{v}_{2}(x^{i},0)-v_{20}^i|^2+...),				
		\end{aligned}
	\end{split}
\end{equation}
\begin{equation}
	\begin{split}
		\begin{aligned}
			loss_{b}=&\frac{1}{N_{b}}\sum_{j=1}^{N_{b}}(|\hat{u}_{1}(x^{j},t^{j})-u_{1b}^j|^2+|\hat{v}_{1}(x^{i},t^{j})-v_{1b}^j|^2\\
			&+|\hat{u}_{2}(x^{j},t^{j})-u_{2b}^j|^2+|\hat{v}_{2}(x^{j},t^{j})-v_{2b}^j|^2+...),			
		\end{aligned}
	\end{split}
\end{equation}
\begin{equation}
	\begin{split}
		\begin{aligned}
			loss_{f}=&\frac{1}{N_{f}}\sum_{p=1}^{N_{f}}(|f_{1u}(x_{f}^p,t_{f}^p)|^2+|f_{1v}(x_{f}^p,t_{f}^p)|^2+|(f_{2u}(x_{f}^p,t_{f}^p)|^2+|f_{2v}(x_{f}^p,t_{f}^p)|^2+...),				
		\end{aligned}
	\end{split}
\end{equation}
here $ \{(x^{i},0),u_{10}^i,v_{10}^i,u_{20}^i,v_{20}^{i},...\}_{i=1}^{N_{0}} $ denotes the initial training data  $ \tau_{0} $ and the boundary training data $ \tau_{b} $ is $ \{(x^{j},t^{j}),u_{1b}^j,v_{1b}^j,u_{2b}^j,v_{2b}^{j},...\}_{j=1}^{N_{b}} $ of Eq.~\eqref{eq2}. What's more, $ \{(x_{f}^p,t_{f}^p)\}_{p=1}^{N_{f}} $ is the collocation points $ \tau_{f} $ in space and time for the residuals $ f_{1u}(x_{f}^p,t_{f}^p),f_{1v}(x_{f}^p,t_{f}^p),...,f_{nv}(x_{f}^p,t_{f}^p)$. $ N_{0},N_{b} $ and $ N_{f} $ are the numbers of auxiliary coordinates required to calculate the $ Loss $ terms in the initial conditions, boundary conditions and governing equations, respectively. $\hat{u}_{1},\hat{v}_{1},..., \hat{v}_{n} $ are the predicted output data through multi-layer neural networks.

\begin{figure}[ht]
\begin{center}
\includegraphics[scale=0.2]{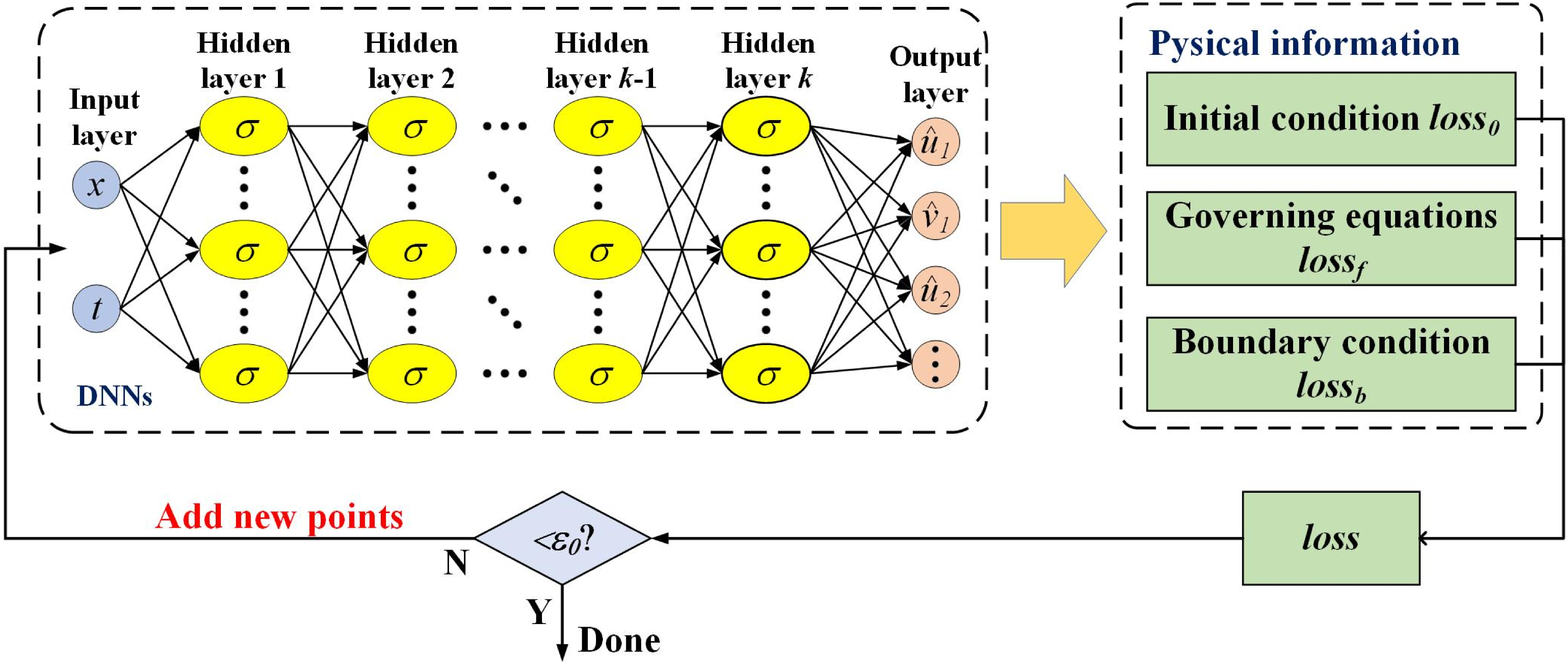}
\caption{Schematic of RAR-PINN for solving the $N$-coupled nonlinear equation.}
\label{figur1}
\end{center}
\end{figure}

As we discussed above, the collocation points for residuals are usually randomly distributed in the solution domain. Furthermore, all randomly sampled point locations are generated using a space filling Latin Hypercube Sampling strategy \cite{1987Large}. This sampled method works well  for many equations, however, it may not be efficient for some PDEs that exhibit solutions with steep gradients \cite{lu2021deepxde}. For better results, we use a RAR-PINN algorithm to improve the distribution of residual points during process. The modified algorithm structure for solving Eq.~\eqref{eq2} is shown in Figure \ref{figur1}, which also applies to other coupled nonlinear equations.

We can divide this structure into three parts, including deep neural networks (DNNs), physical information and loss function. The part of DNNs is the traditional fully connected feedforward neural networks. In detail, the input layer consists of two neurons $ x $ and $ t $, the number of hidden layers is $ k $, and the output layer consists of $ 2n $ neurons $\hat{u}_{1},\hat{v}_{1},\hat{u}_{2},...,\hat{v}_{n}$. We can compute the derivatives of the neural network outputs with respect to the neural network inputs time $ t $ and space $ x $ by AD, and incorporate them into residual as a part of the loss function. AD, also known as algorithmic differentiation or simply ``autodiff'', is a technique similar to but more general than backpropagation for efficiently and accurately evaluating the derivative of a function expressed as a computer program \cite{baydin2018automatic}. TensorFlow \cite{abadi2016tensorflow} or PyTorch \cite{paszke2017automatic} is used to implement AD of the parameters in a neural network.
In the physical information section, $ loss_{0} $ is determined by the initial condition, the governing equation controls the numerical value of $ loss_{f} $, and the numerical value of $ loss_{b} $ is determined by the boundary condition. In addition, the value of the loss function can be obtained from Eq.~\eqref{Loss}. Furthermore, in order to improve the efficiency of the neural networks, we judge whether the mean residual of the PDE is less than a threshold $\varepsilon _{0} $. If it is less than $\varepsilon _{0} $, the training ends. Otherwise, we add new points to continue the training until the mean residual is less than $ \varepsilon_{0} $, as shown in algorithm \ref{1}. RAR-PINN method adds more degrees of freedom when the solution PDE residual is high. \\ \hspace*{\fill} \\
\begin{algorithm}[H]
	\caption{RAR-PINN for solving Eq.~\eqref{eq2}}
	\KwIn{Initial training data set $ \tau_{0}, \tau_{b}$ and $ \tau_{f} $, number of new points added $ m $, number of iterations $ n $ and a tolerance $ \varepsilon_{0} $; }

	\While{$ err < \varepsilon_{0} $}{
		
	   -- Calculate the mean coupled nonlinear differential partial equation residual $ err $, specifically,\\
	   \centering	
    	$err=\dfrac{1}{N}(|f_{1u}|+|f_{1v}|+|f_{2u}|+|f_{2v}|+...);$

       -- Select the $ m $ points with the largest residual\;

       -- Add the $ m $ new points to the training data set\;
	}
	\KwOut{Predictive solutions of Eq.~\eqref{eq1}.}
     \label{1}
\end{algorithm}

\section{Data-driven vector $N$-soliton solutions of the CGNLS equation }
\subsection{Vector one-soliton solutions}

In the following, we will take Eq.~\eqref{equation} as an example. Vector one-soliton solutions of  Eq.~\eqref{equation} was derived in \cite{2013Soliton} via the Hirota bilinear method. However, the data-driven vector one-soliton solutions of Eq.~\eqref{equation} have not been investigated so far. In order to assess the accuracy of our method, here, we give the exact vector one-soliton solution for System \eqref{equation} as \cite{2013Soliton}
\begin{equation}
	\left\{
	\begin{aligned}
	&h_{1}=\frac{ae^{\theta}}{1+e^{ \theta + \theta ^{*} +c}}=\frac{a}{2} e^{-\frac{c}{2}}e^{i\theta_{I}}{\rm sech}(\theta_{R}+\frac{c}{2}),\\
	&h_{2}=\frac{be^{\theta}}{1+e^{\theta+\theta^{*}+c}}=\dfrac{b}{2} e^{-\frac{c}{2}}e^{i\theta_{I}}{\rm sech}(\theta_{R}+\frac{c}{2}),\label{one}
	\end{aligned}
	\right.
\end{equation}
with $ \theta=kx+ik^{2}t $ and $ e^{c}= \frac{|a|^{2}\alpha+|b|^{2}\beta+a_{1}b^{*}\gamma+a^{*}b\gamma^{*}}{(k+k^{*})^2}$, where $ \frac{c}{2} $ represents the localization position of
one-solitons; $ a $, $ b $ and $ k $ are complex constants. Here and elsewhere, the subscripts $ R $ and $ I $ represent the real and imaginary parts, respectively.

\begin{figure}[htbp]
	\centering
		\includegraphics[scale=0.4]{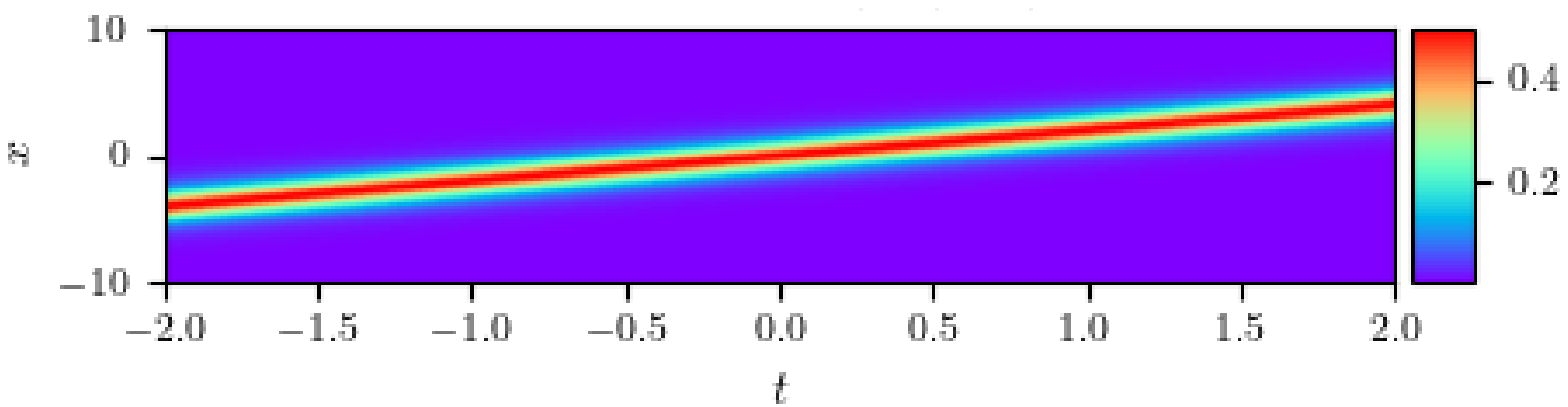}\hspace{1cm}
		\includegraphics[scale=0.4]{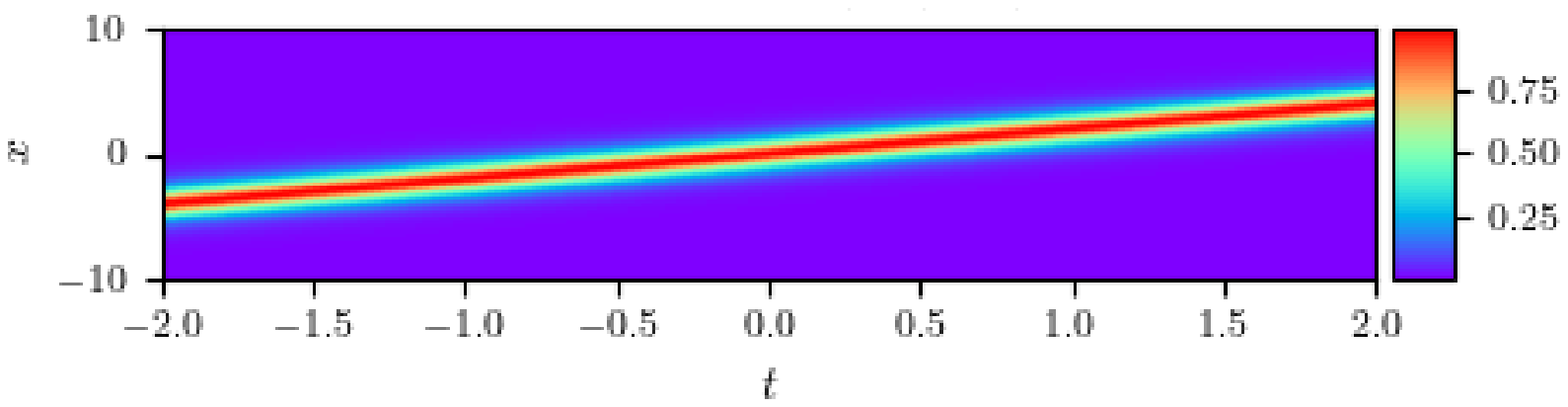}
		{\footnotesize\hspace{4cm} (a) Exact solution $ |h_{1}(x,t)| $ \hspace{4.5cm}(b) Exact solution $ |h_{2}(x,t)| $ } \\
		\vspace{5mm}
		\includegraphics[scale=0.4]{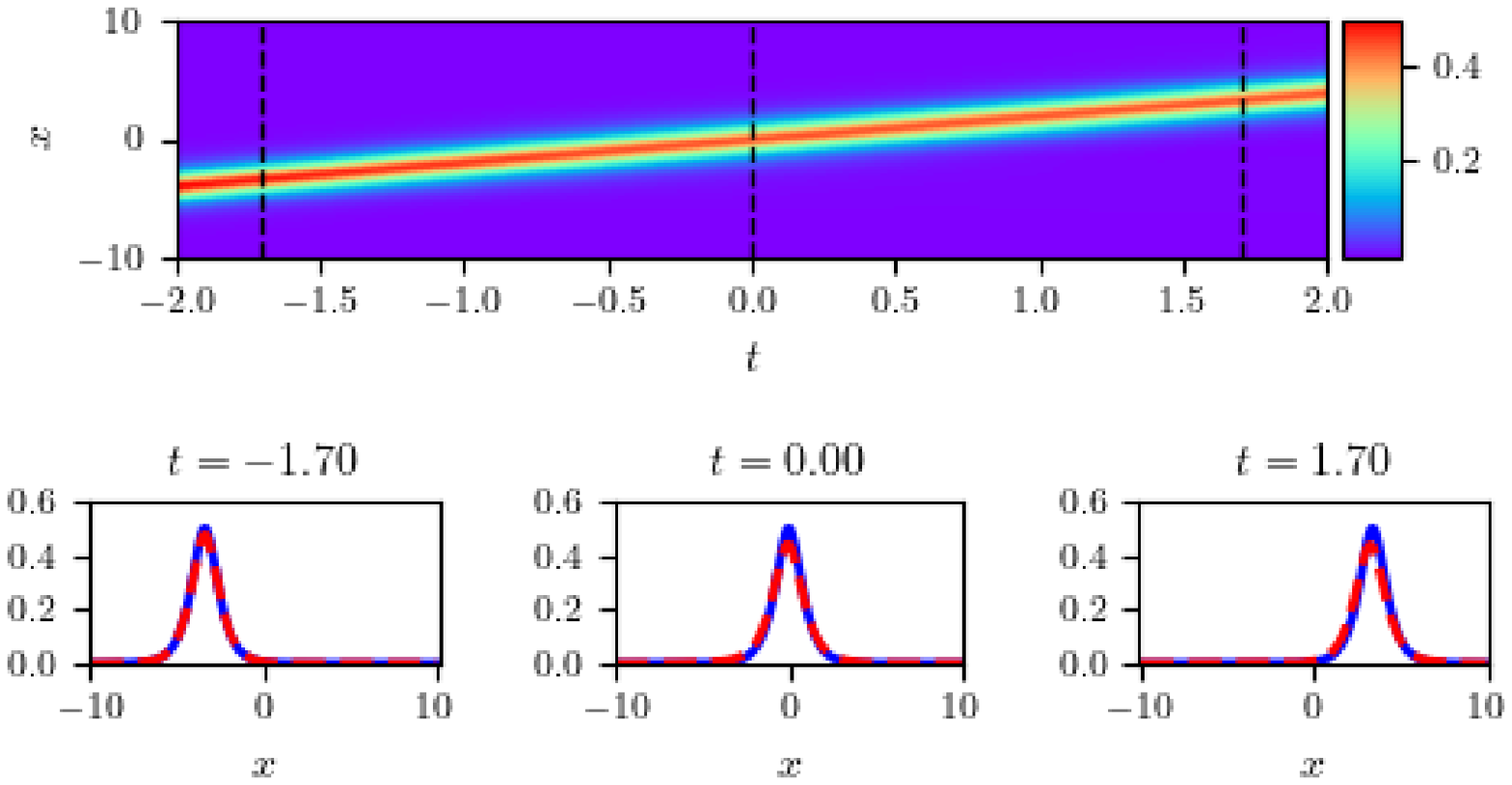}\hspace{1cm}
		\includegraphics[scale=0.4]{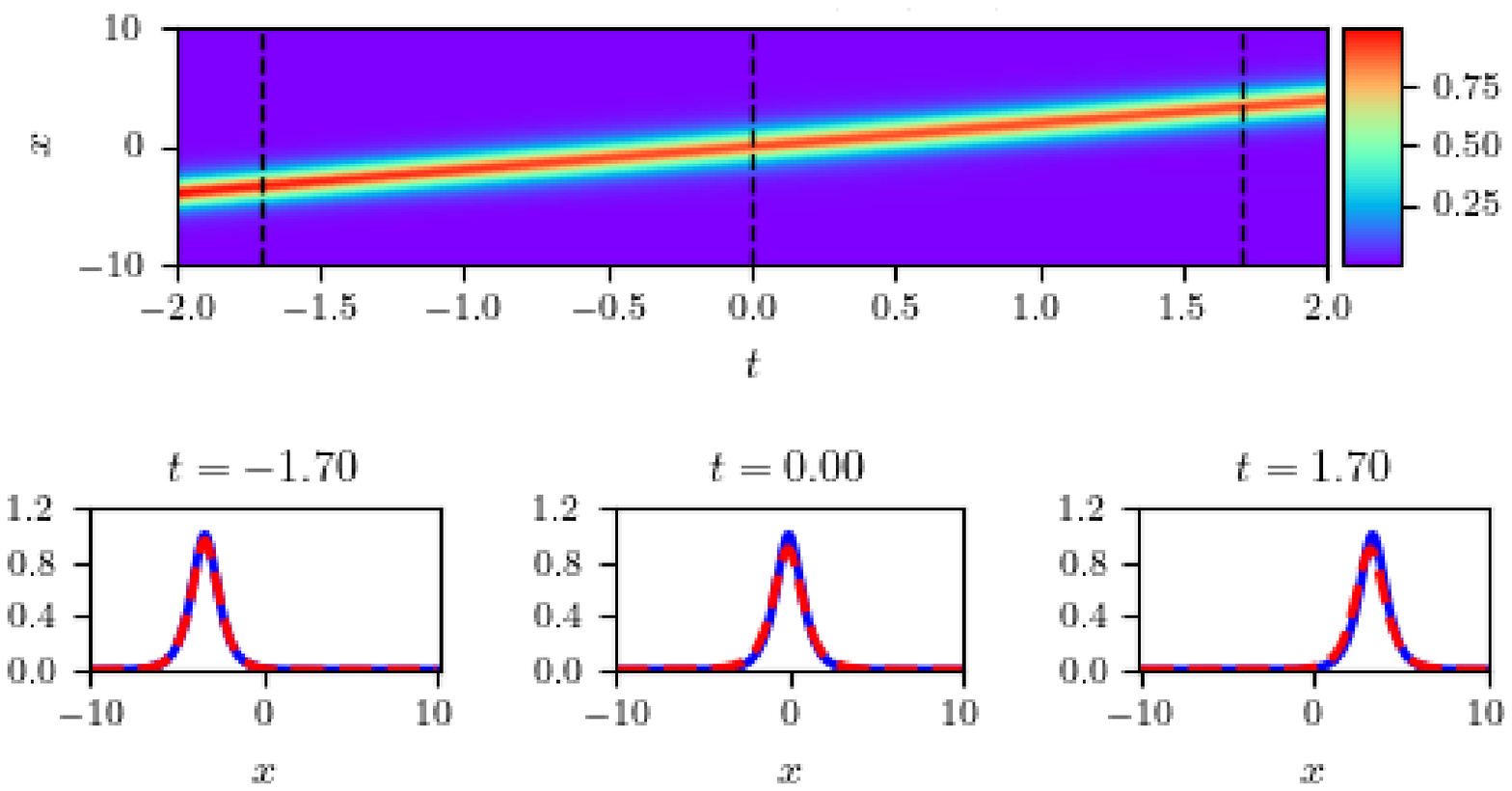}
		{\footnotesize\hspace{3cm} (c) Trational PINN solution $ |h_{1}(x,t)| $ \hspace{3cm}(d) Trational PINN solution $ |h_{2}(x,t)| $ } \\
		\vspace{5mm}
		\includegraphics[scale=0.4]{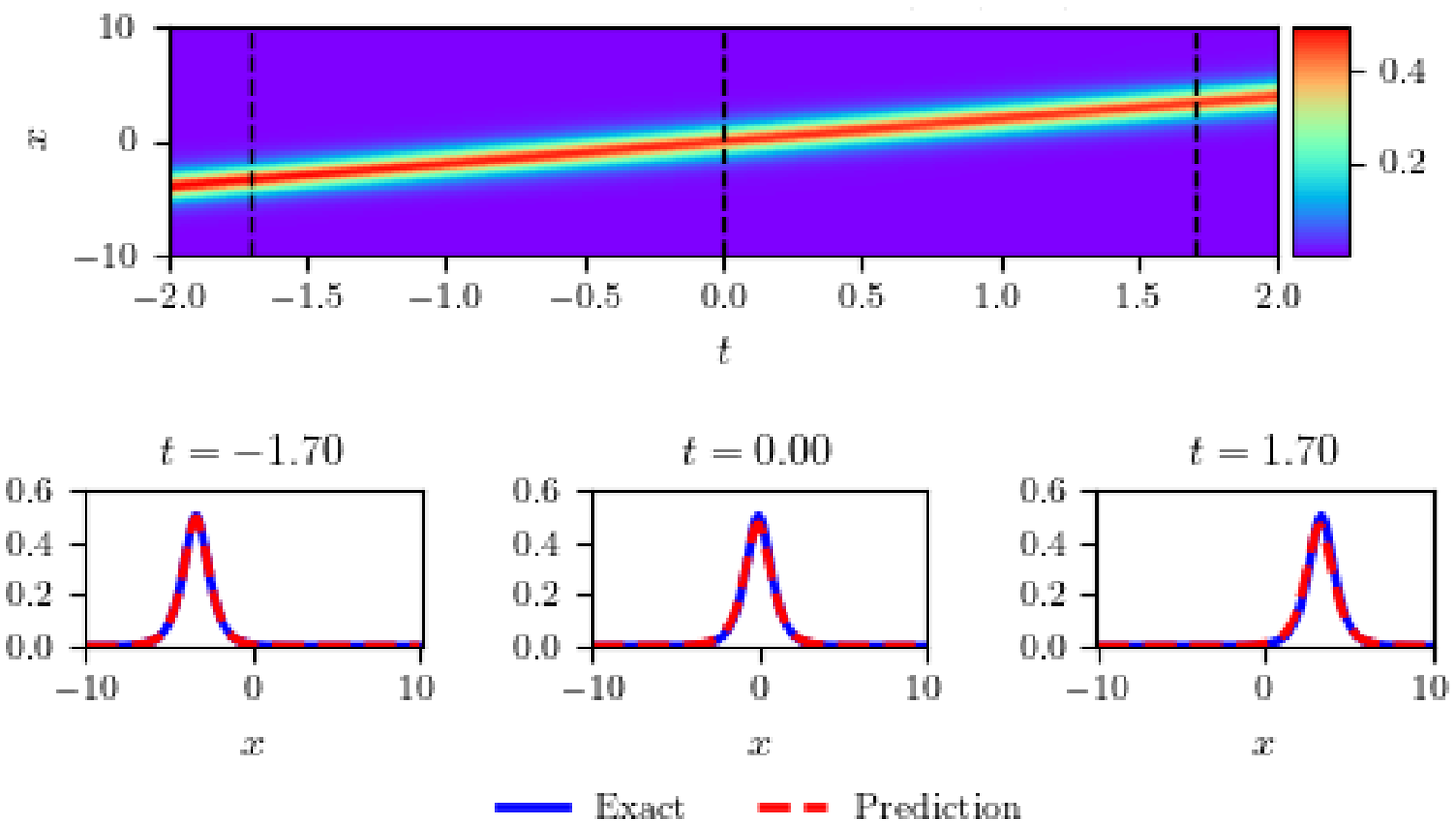}\hspace{1cm}
		\includegraphics[scale=0.4]{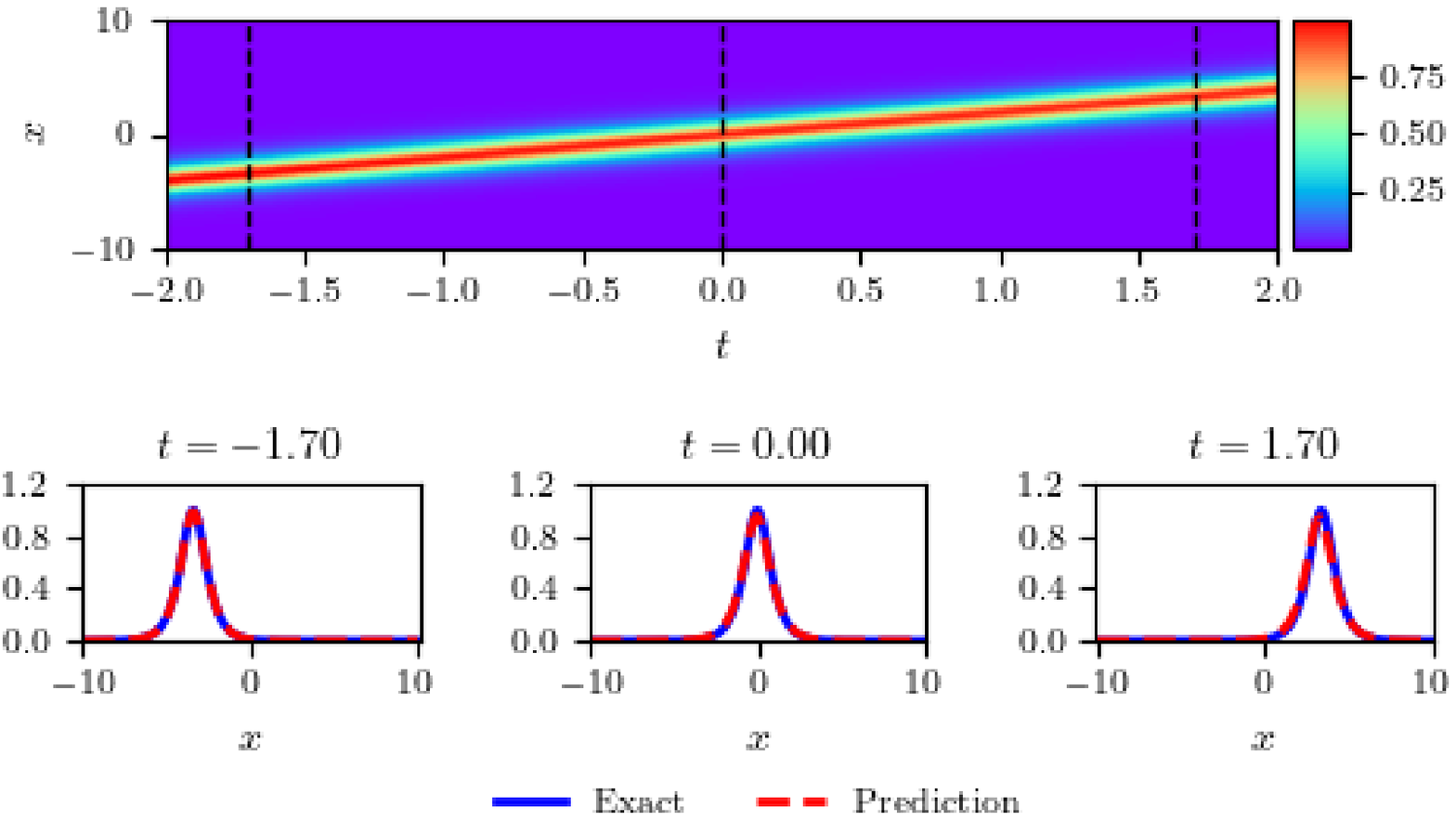}
		{\footnotesize\hspace{3.5cm} (e) RAR-PINN solution $ |h_{1}(x,t)| $ \hspace{3.5cm}(f) RAR-PINN solution $ |h_{2}(x,t)| $ }  \\
		\caption{The vector one-soliton solution $ h_{1}(x,t) $ and $ h_{2}(x,t) $: Comparison among the exact solutions (panels (a) and (b)), predictive  solutions via TPINN (panels (c) and (d)) and predictive results by RAR-PINN (panels (e) and (f)) at different propagation times $ t=-1.70, 0.0 $ and 1.70. The left and right graphs represent the density diagrams and the profiles at different times for the one-soliton solutions $ h_{1}(x,t) $ and $ h_{2}(x,t) $, respectively.
		}
		\label{figure2}
\end{figure}

As an example, by choosing $ a=1 $, $ b=2 $, $ k=1.5+i $ and $ \alpha=\beta=\gamma=1 $ in Solution \eqref{one}
the corresponding initial condition is obtained as
\begin{equation}
	\left\{
	\begin{aligned}
		&h_{10}(x)=h_{1}(x,-2)=\frac{e^{(1.5+i)x+6-2.5i}}{1+e^{3x+12}} =\frac{1}{2}e^{(1.5i-1)x+6i+2.5}{\rm sech}(1.5x+6),\\
		&h_{20}(x)=h_{2}(x,-2)=\frac{2e^{(1.5+i)x+6-2.5i}}{1+e^{3x+12}}= e^{(1.5i-1)x+6i+2.5}{\rm sech}(1.5x+6),
	\end{aligned}
	\right.
\end{equation}
and the Dirichlet-Neumann periodic boundary condition is expressed as
\begin{equation}
	\left\{
	\begin{aligned}
		&h_{1}(-10,t)=h_{1}(10,t),h_{2}(-10,t)=h_{2}(10,t),t\in[-2,2], \\
		&h_{1x}(-10,t)=h_{1x}(10,t),h_{2x}(-10,t)=h_{2x}(10,t).
	\end{aligned}
	\right.
\end{equation}
The training data set can be obtained by using random sampling. More specifically, the computation domain $ [-10,10]\times[-2,2] $ is divided into $ [300\times201] $ data points, vector one-soliton solutions $ h_{1} $ and $ h_{2} $ are discretized into 201 snapshots on regular space-time grid with $ \bigtriangleup t=0.02 $. A small training data set containing initial data is generated by randomly extracting $ N_{0}=50 $ from original data set. Similarly, the number of boundary points is $ N_{b}=50 $. Initially, 4000 points are selected as the residual points, and then we add 10 more new points adaptively via RAR with $ m=5 $ and $ \varepsilon_{0}=0.01 $. To compare the RAR-PINN solutions with TPINN solutions, we obtain the TPINN solutions based on $ N_{f}=4010 $ points by the Latin Hypercube Sampling method \cite{1987Large}. After giving a data set of initial points, boundary points and residual points, the latent vector one-solitons $ h_{1}(x,t) $ and $ h_{2}(x,t) $ can be successfully learned by minimizing the loss function Eq.~\eqref{Loss}. In order to minimiz the loss function, we use 6 hidden layers of 32 neurons in each layer to numerically simulate soliton solutions. The results are displayed in Figure \ref{figure2}, which demonstrates the effectiveness of the RAR-PINN method.

From Figure \ref{figure2}, we can clearly observe the reconstructed vector one-soliton space-time dynamics. The density plots with the corresponding peak scale for diverse dynamics include exact dynamics, learned dynamics with TPINN and learned dynamics with RAR-PINN, respectively. Specifically, Figure \ref{figure2}(a) and Figure \ref{figure2}(b) indicate the magnitude of the exact vector one-soliton solutions $|h_{1}(x,t)|=\sqrt{u_{1}^{2}(x,t)+v_{1}^{2}(x,t)}$ and $|h_{2}(x,t)|=\sqrt{u_{2}^{2}(x,t)+v_{2}^{2}(x,t)}$; Figure \ref{figure2}(c) and Figure \ref{figure2}(d) show the the predicted vector one-soliton solutions via TPINN; and Figure \ref{figure2}(e) and Figure \ref{figure2}(f) demonstrate the results via RAR-PINN. To make a more detailed comparison, the relative $ \mathbb{L}_{2} $ errors of the TPINN model are $9.189407\times10^{-2}$ for $h_{1}(x,t)$ and $8.453297\times10^{-2}$ for $h_{2}(x,t)$; while the relative $\mathbb{L}_{2}$ errors of the RAR-PINN model are $4.979434\times10^{-2}$ for $h_{1}(x,t)$ and $ 4.661114\times10^{-2}$ for $h_{2}(x,t)$. Besides, Figure \ref{figure2} demonstrates
the comparison of the exact vector one-soliton solution (blue curve) and the predicted one (red curve) at different times.

\begin{figure}[htbp]
	\centering
    	\includegraphics[scale=0.06]{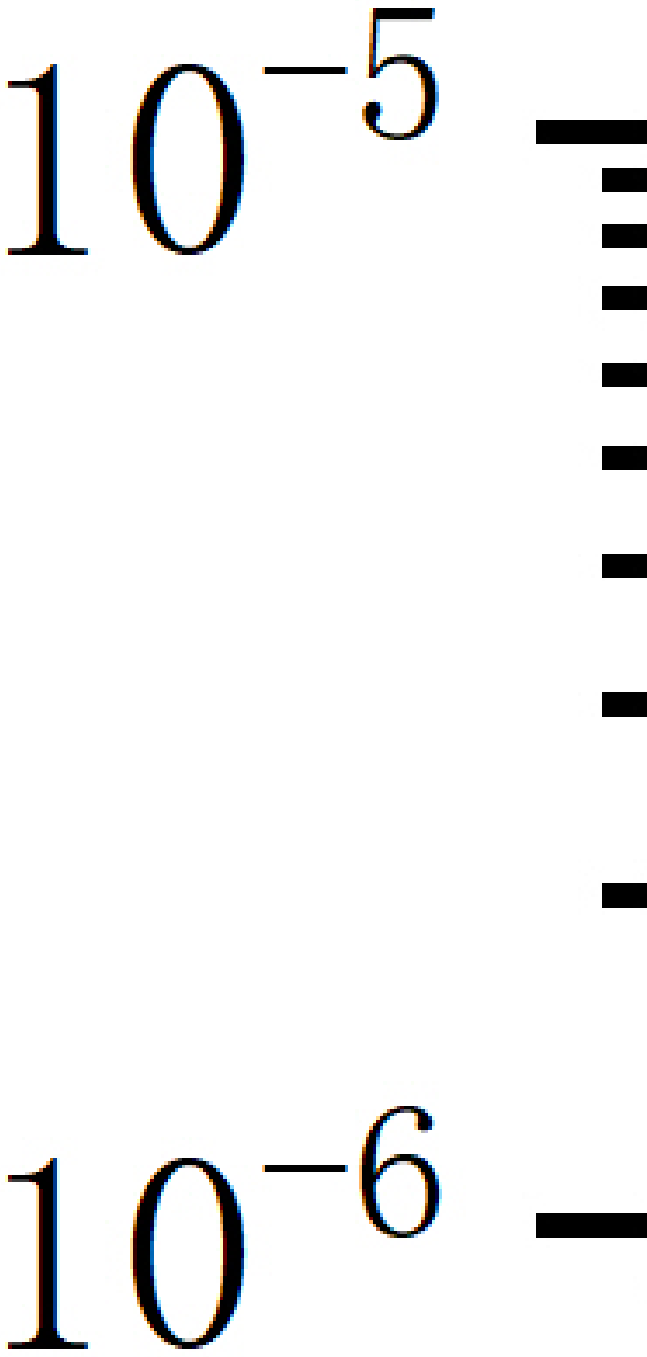}\hspace{1cm}
    	\includegraphics[scale=0.06]{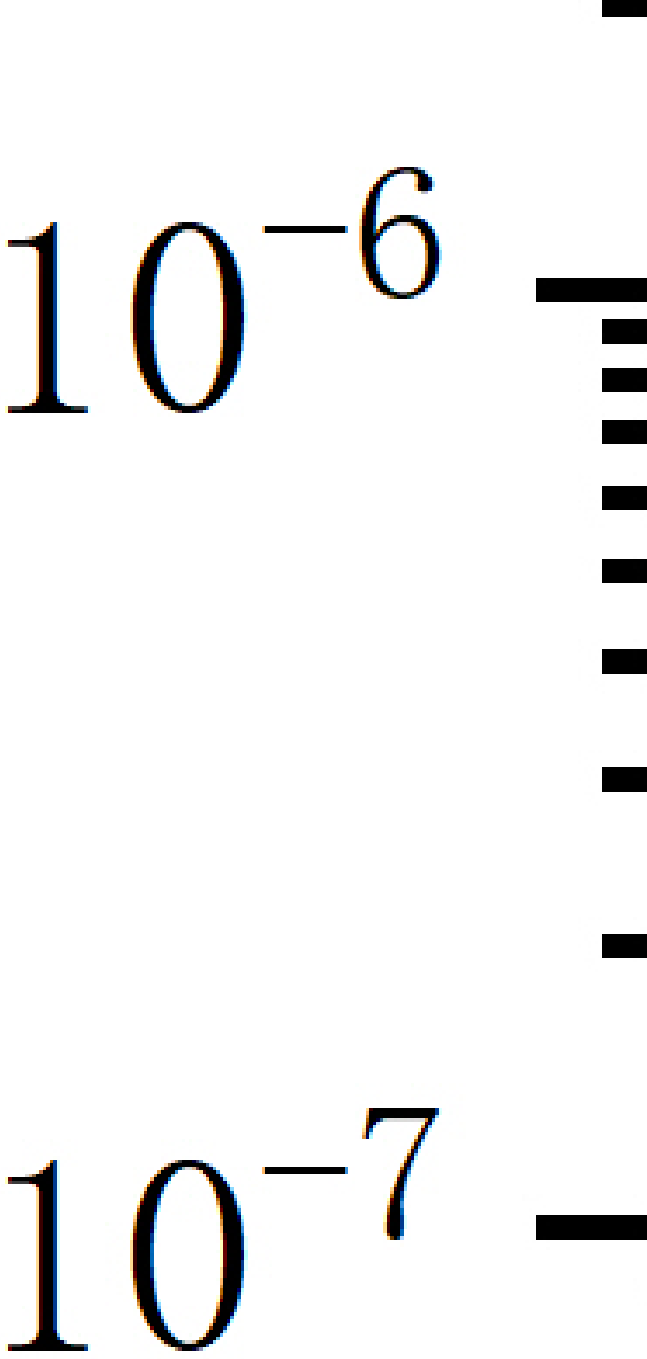}
    	{\footnotesize\hspace{9cm} (a) \hspace{8cm}(b) } \\
    	\includegraphics[scale=0.06]{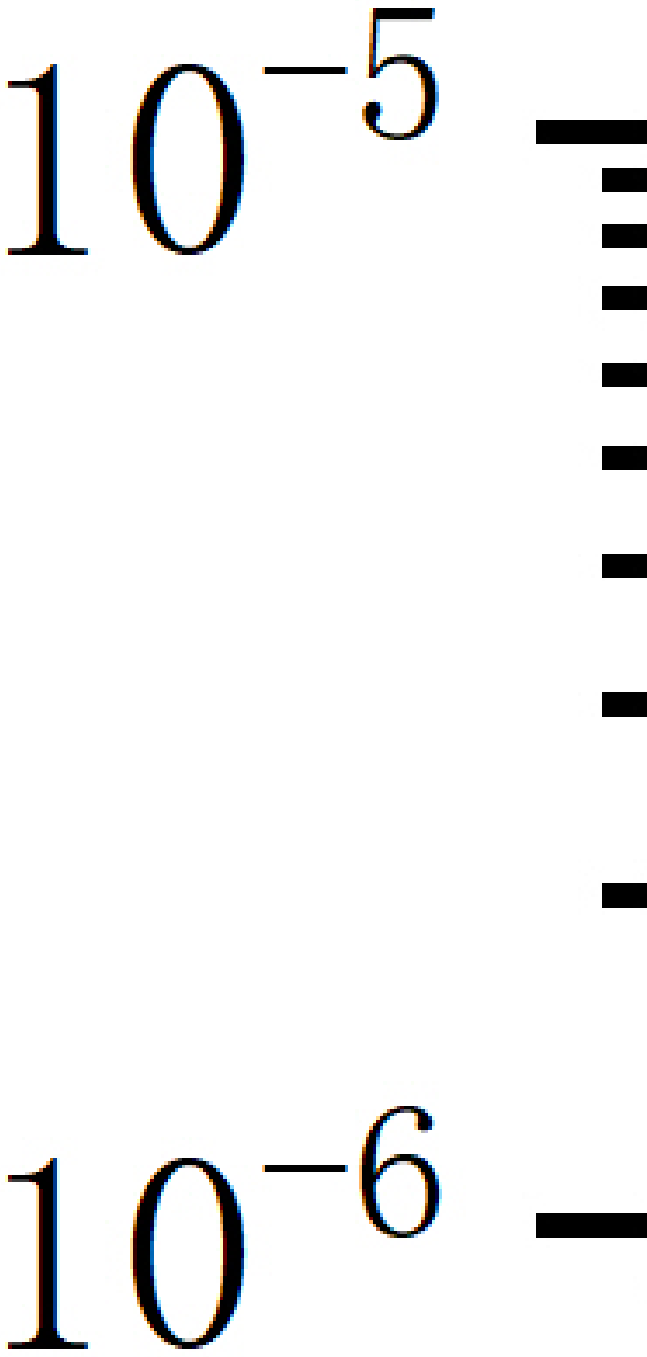}\hspace{1cm}
    	\includegraphics[scale=0.06]{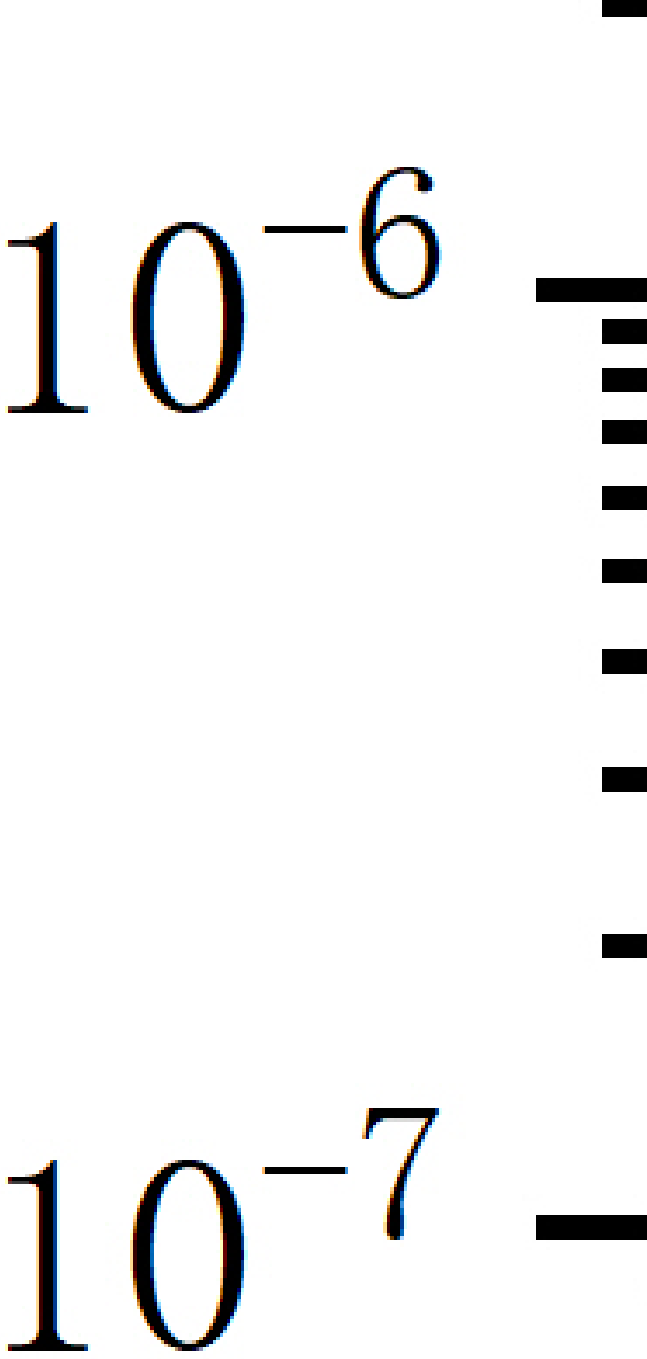}
    	{\footnotesize\hspace{9cm} (c) \hspace{8cm}(d) } \\  	
		\caption {The loss function curve figures of one-soliton solutions $ h_{1}(x,t) $ and $ h_{2}(x,t) $ arising from (a) TPINN algorithm with the 10000 iterations Adam optimization; (b) TPINN algorithm with the 9864 iterations L-BFGS optimization; similarly, (c) RAR-PINN algorithm with the 10000 iterations Adam optimization; (d) RAR-PINN algorithm with the 13506 iterations L-BFGS optimization.}
		\label{figure3}
\end{figure}

Specifically, to obtain the optimal parameters to minimize the loss function,
we first use the Adam optimizer \cite{diederik2014adam} for 10,000 iterations, followed by the L-BFGS optimizer \cite{1989On}. Importantly, the L-BFGS method has global convergence on uniformly convex problems and it is very competitive due to its low iteration cost \cite{1989On}. The specific loss function curve is shown in Figure \ref{figure3}. Figure \ref{figure3}(a) and (b) represent the loss function curves via TPINN, in detail, Figure \ref{figure3}(a) shows the loss function descent curve for 10,000 iterations using the Adam optimizer and Figure \ref{figure3}(b) demonstrates the loss function curve using the L-BFGS optimizer. From Figure \ref{figure3}(b) the initial condition error is $loss_{0}=2.2284491\times10^{-6}$, the boundary condition error is $ loss_{b}=4.350938\times10^{-7}$, the collocation points error is $loss_{f}=1.6587395\times10^{-5}$, and the overall error is $ Loss= 1.9250938\times10^{-5}$. Figure \ref{figure3}(c) and Figure \ref{figure3}(d) illustrate the loss function curves with RAR-PINN, specifically, Figure \ref{figure3}(c) indicates the loss function curve for 10,000 iterations using the Adam optimization and Figure \ref{figure3}(d) shows the loss function curve using the L-BFGS optimization. From Figure \ref{figure3}(d) the initial condition error is $loss_{0}=6.016333\times10^{-7}$, the boundary condition error is $ loss_{b}=1.3800212\times10^{-7}$, the collocation points error is $loss_{f}=7.81987\times10^{-6}$ and the overall error is $Loss= 8.559506\times10^{-6}$. By comparing the loss functions of the two methods, the RAR-PINN method is more efficient for solving the vector one-soliton solutions of Eq.~\eqref{equation}.

\begin{figure}[htbp]
	\centering
	\includegraphics[scale=0.06]{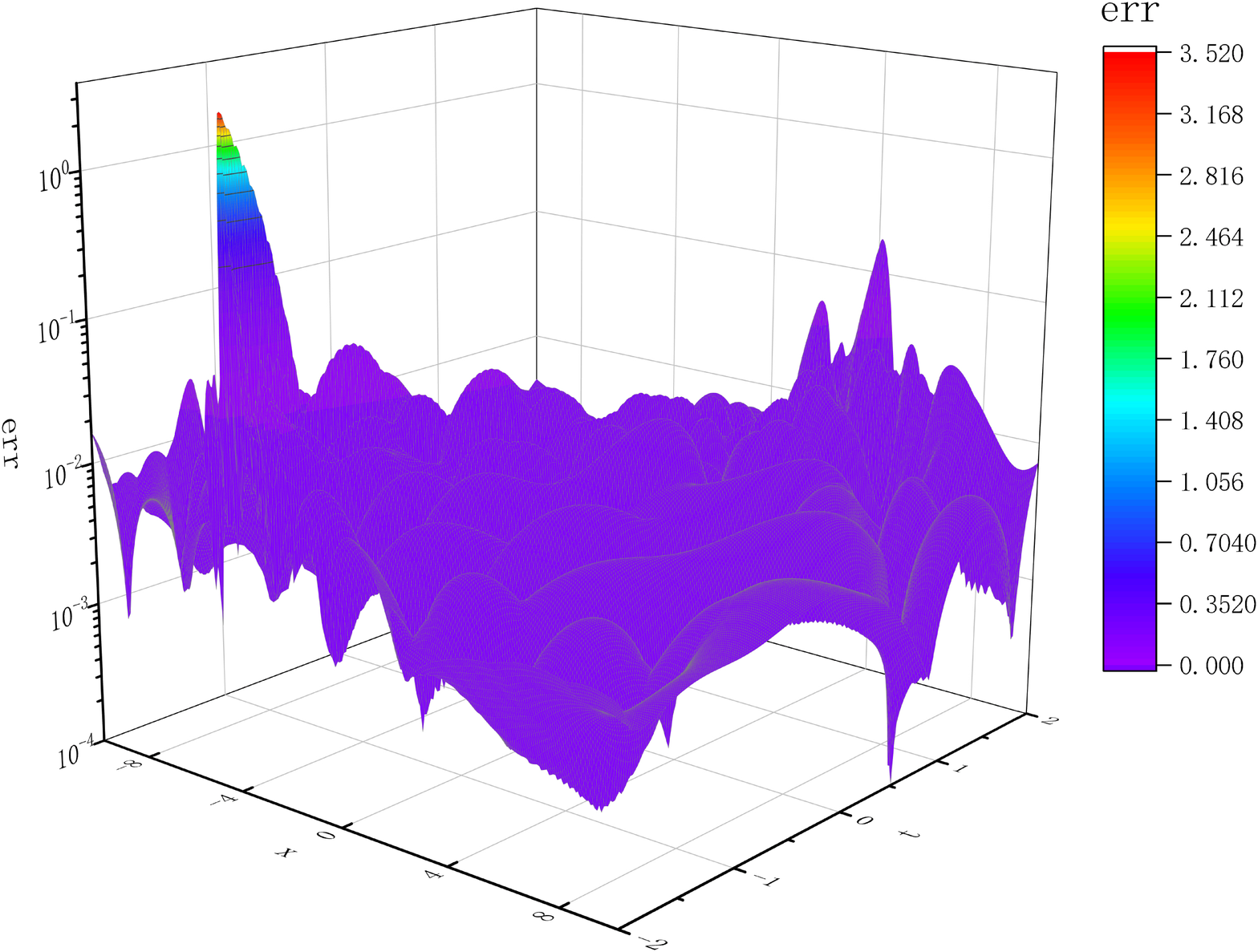}\hspace{1cm}
	\includegraphics[scale=0.06]{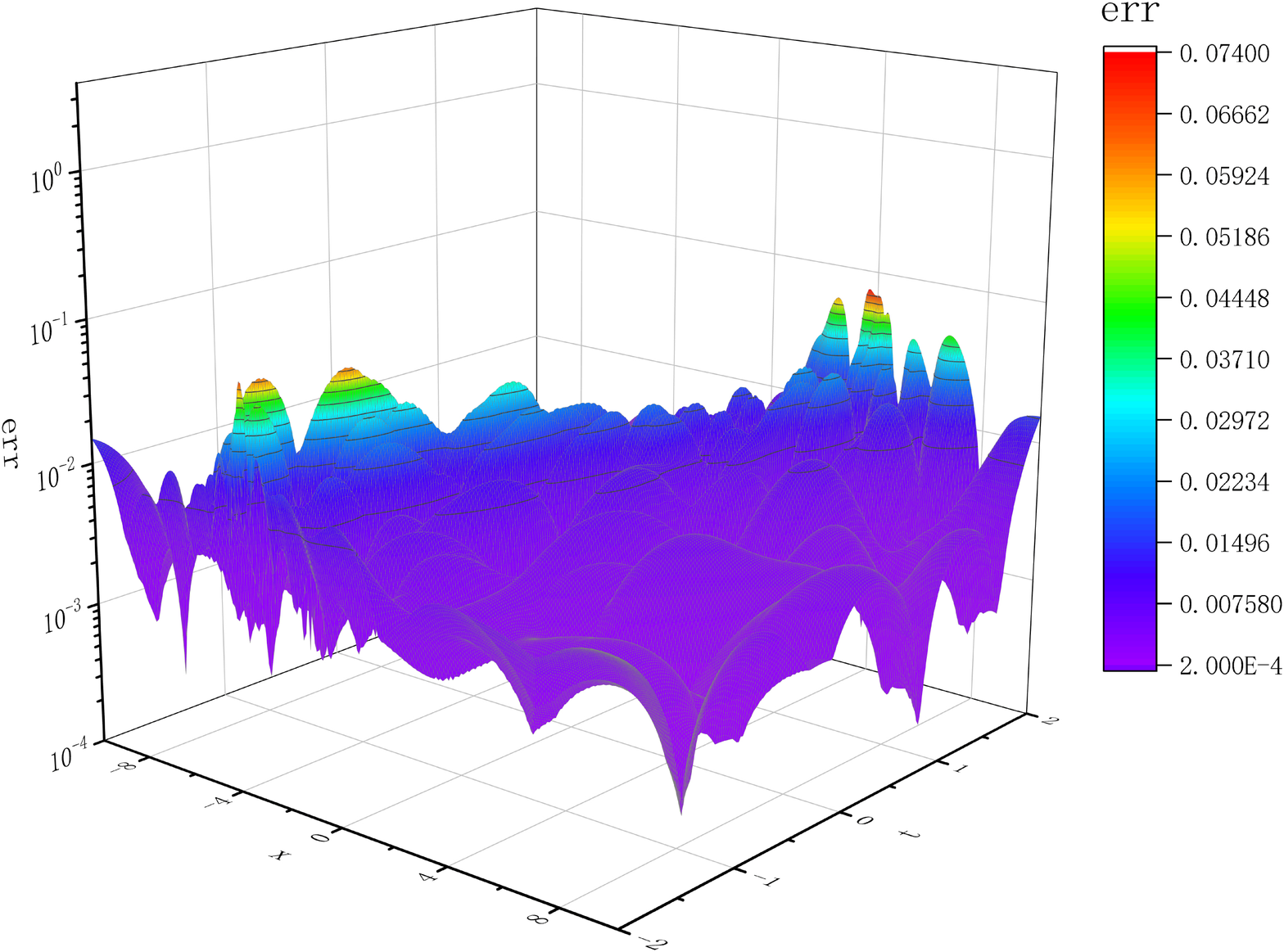}
	{\footnotesize\hspace{9cm} (a) \hspace{8cm}(b) } 	
	\caption{Three-dimensional stereogram of residual error: (a) three-dimensional stereogram residual error without adding adaptive points; (b) three-dimensional stereogram residual error after adding adaptive points.
}
	\label{figure4}
\end{figure}

The three-dimensional diagram of the residual error is shown in Figure \ref{figure4}. From Figure \ref{figure4}(a), we can see that the residuals of some positions before the adaptive pointing are relatively large up to 3.520, and from Figure \ref{figure4}(b), the residuals of all positions after the adaptive pointing can be controlled within 0.08. Combining the above numerical simulation results, we can find that RAR-PINN is more accurate than TPNN for solving the vector one-soliton solutions of Eq.~\eqref{equation}.

\subsection{Vector two-soliton solutions}

Through the Hirota bilinear method, the vector two-soliton solutions can be derived via symbolic computation. It is found that the vector two-soliton solutions admit the inelastic and elastic collisions depending on the choice of solitonic parameters. Simplifying the resultant expressions we have the exact vector two-soliton solutions of the form~\cite{VISHNUPRIYA2016366}

\begin{equation}
	\left\{
	\begin{aligned}
		&h_{1}=\frac{\xi_{1}^{(1)}e^{\eta_{1}}+\xi_{2}^{(1)}e^{\eta_{2}}+e^{\eta_{1}+\eta_{1}^{*}+\eta_{2}+\delta_{11}}+e^{\eta_{1}+\eta_{2}^{*}+\eta{2}+\delta_{21}}}{D},\\
		&h_{2}=\frac{\xi_{1}^{(2)}e^{\eta_{1}}+\xi_{2}^{(2)}e^{\eta_{2}}+e^{\eta_{1}+\eta_{1}^{*}+\eta_{2}+\delta_{12}}+e^{\eta_{1}+\eta_{2}^{*}+\eta{2}+\delta_{22}}}{D},\label{two}
	\end{aligned}
	\right.
\end{equation}
and
\begin{equation}
	\begin{split}
		\begin{aligned}
		D=1+e^{\eta_{1}+\eta_{1}^{*}+R_{1}}+e^{\eta_{1}+\eta_{2}^{*}+\delta_{0}}+e^{\eta_{1}^{*}+\eta_{2}+\delta_{0}^{*}}+e^{\eta_{2}+\eta_{2}^{*}+R_{2}}+e^{\eta_{1}+\eta_{1}^{*}+\eta_{2}+\eta_{2}^{*}+R_{3}},		
		\end{aligned}
\end{split}
\end{equation}
where we have defined the constants as
\begin{equation}
	\begin{split}	
	&\eta_{1}=k_{1}(x+ik_{1}t), \eta_{2}=k_{2}(x+ik_{2}t), e^{\delta_{0}}=\frac{\varphi_{12}}{k_{1}+k_{2}^{*}},e^{\delta_{R_{1}}}=\frac{\varphi_{11}}{k_{1}+k_{1}^{*}},\\
	&e^{\delta_{R_{2}}}=\frac{\varphi_{22}}{k_{2}+k_{2}^{*}},
	e^{\delta_{R_{3}}}=\frac{|k_{1}-k_{2}|^{2}(\varphi_{11}\varphi_{22}-\varphi_{12}\varphi_{21})}{(k_{1}+k_{1}^{*})(k_{2}+k_{2}^{*})|k_{1}+k_{2}^{*}|^{2}},\\
	&e^{\delta_{11}}=\frac{(k_{1}-k_{2})(\xi_{1}^{(1)}\varphi_{21}-\xi_{2}^{(1)}\varphi_{11})}{(k_{1}+k_{1}^{*})(k_{2}+k_{1}^{*})},
	e^{\delta_{21}}=\frac{(k_{2}-k_{1})(\xi_{2}^{(2)}\varphi_{12}-\xi_{1}^{(1)}\varphi_{22})}{(k_{2}+k_{2}^{*})(k_{1}+k_{2}^{*})},\\
	&e^{\delta_{12}}=\frac{(k_{1}-k_{2})(\xi_{1}^{(2)}\varphi_{21}-\xi_{2}^{(2)}\varphi_{11})}{(k_{1}+k_{1}^{*})(k_{2}+k_{1}^{*})},
	e^{\delta_{22}}=\frac{(k_{2}-k_{1})(\xi_{2}^{(2)}\varphi_{12}-\xi_{1}^{(2)}\varphi_{22})}{(k_{2}+k_{2}^{*})(k_{1}+k_{2}^{*})},		
   \end{split}
\end{equation}
and
\begin{equation}
	\varphi_{mn}=\dfrac{\alpha\xi_{m}^{(1)*}+\beta\xi_{m}^{(2)}\xi_{n}^{(2)*}+\gamma\xi_{m}^{(1)}\xi_{n}^{(2)*}+\gamma\xi_{m}^{(1)*}\xi_{n}^{(2)}}{k_{m}+k_{n}^{*}},  m,n=1,2.
	\end{equation}
The two-soliton solutions $ h_{1}(x,t) $ and $ h_{2}(x,t) $ are characterized by arbitrary complex parameters
$ k_{1} $, $ k_{2} $, $ \xi_{1}^{(1)} $, $ \xi_{1}^{(2)} $, $ \xi_{2}^{(1)} $ and $ \xi_{2}^{(2)} $.

\subsubsection{Elastic collision of vector two-soliton solution}

\begin{figure}[htbp]
	\centering
	\includegraphics[scale=0.4]{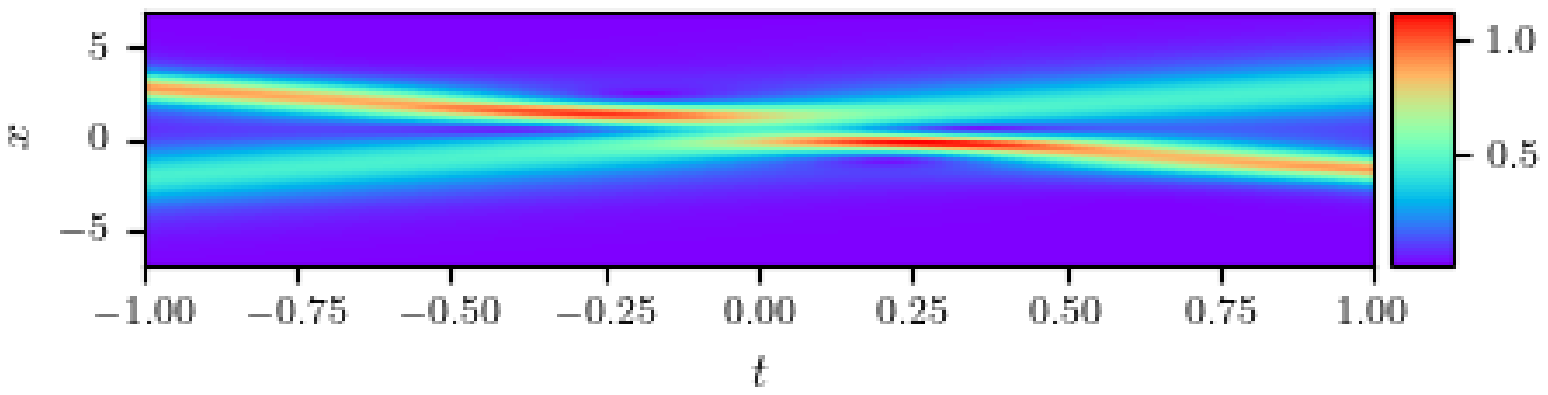}\hspace{1cm}
	\includegraphics[scale=0.4]{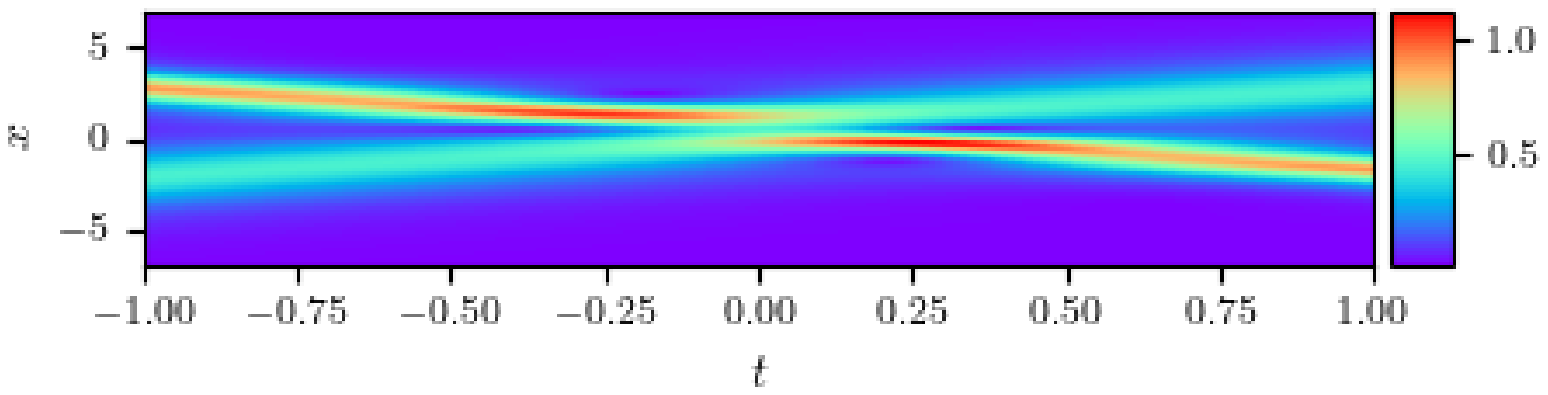}
	{\footnotesize\hspace{4cm} (a) Exact solution $ |h_{1}(x,t)| $ \hspace{4.5cm}(b) Exact solution $ |h_{2}(x,t)| $ } \\
	\vspace{5mm}
	\includegraphics[scale=0.4]{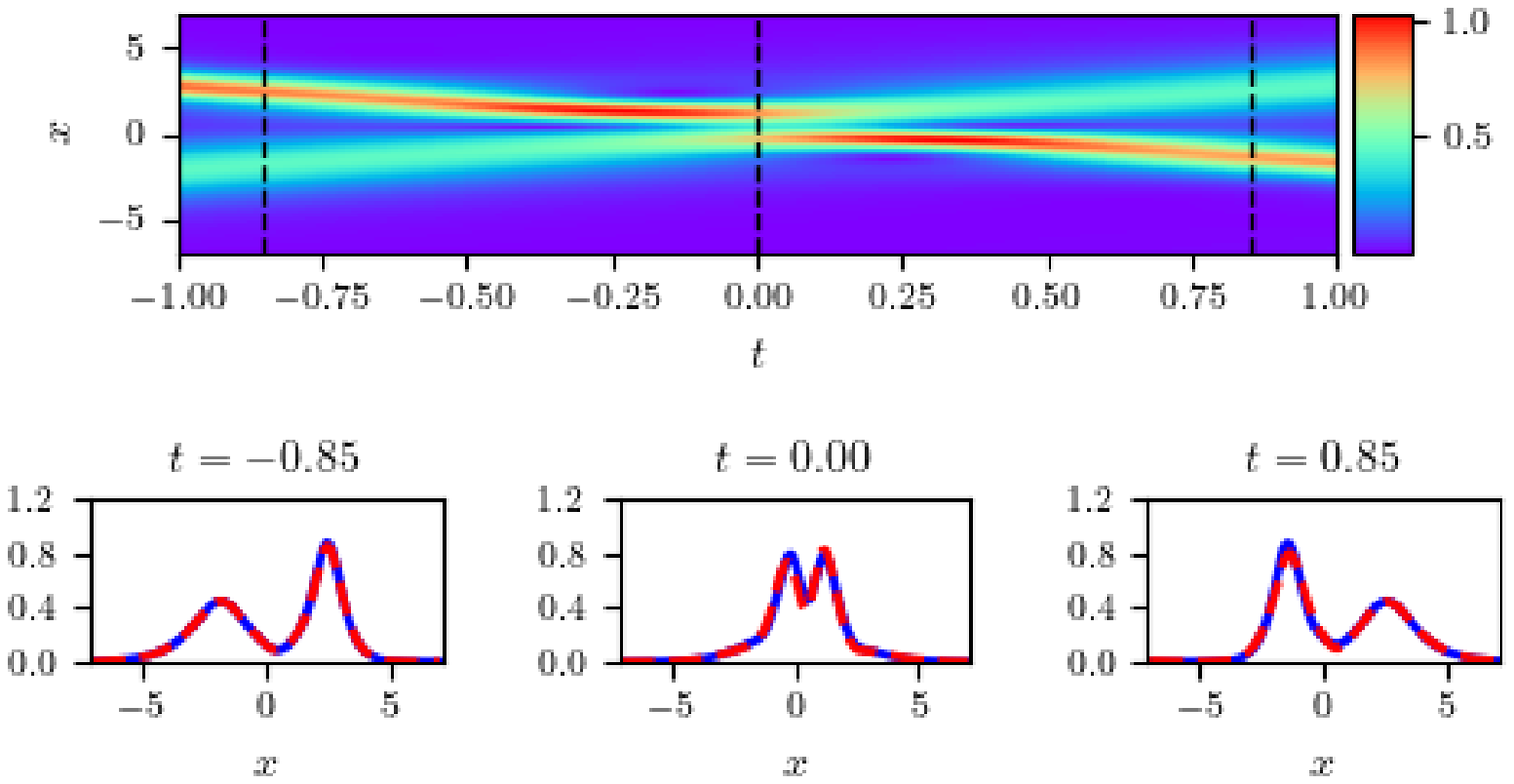}\hspace{1cm}
	\includegraphics[scale=0.4]{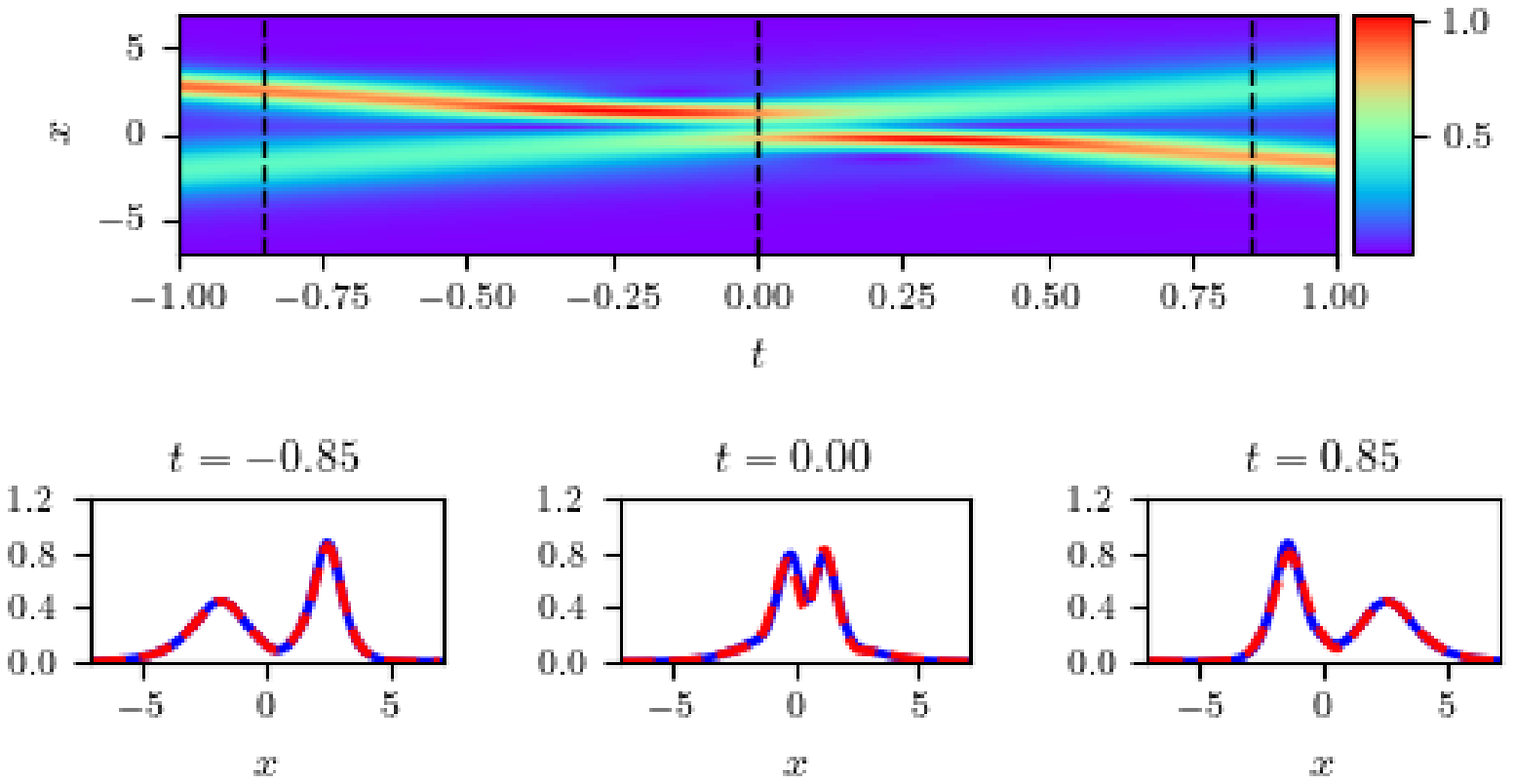}
	{\footnotesize\hspace{3cm} (c) Tractional PINN solution $|h_{1}(x,t)|$ \hspace{3cm}(d) Tractional PINN solution $ |h_{2}(x,t)| $ } \\
	\vspace{5mm}
	\includegraphics[scale=0.4]{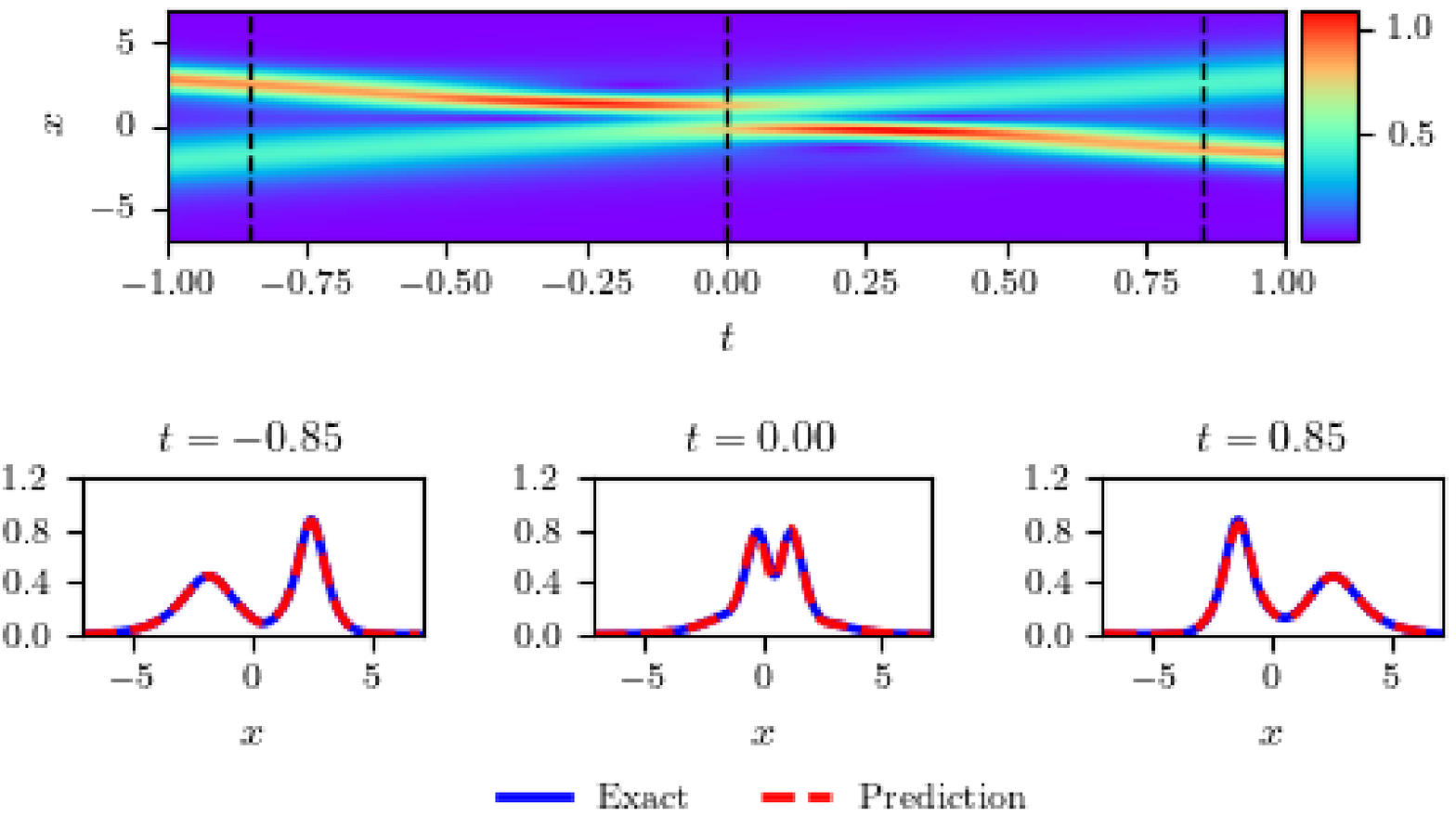}\hspace{1cm}
	\includegraphics[scale=0.4]{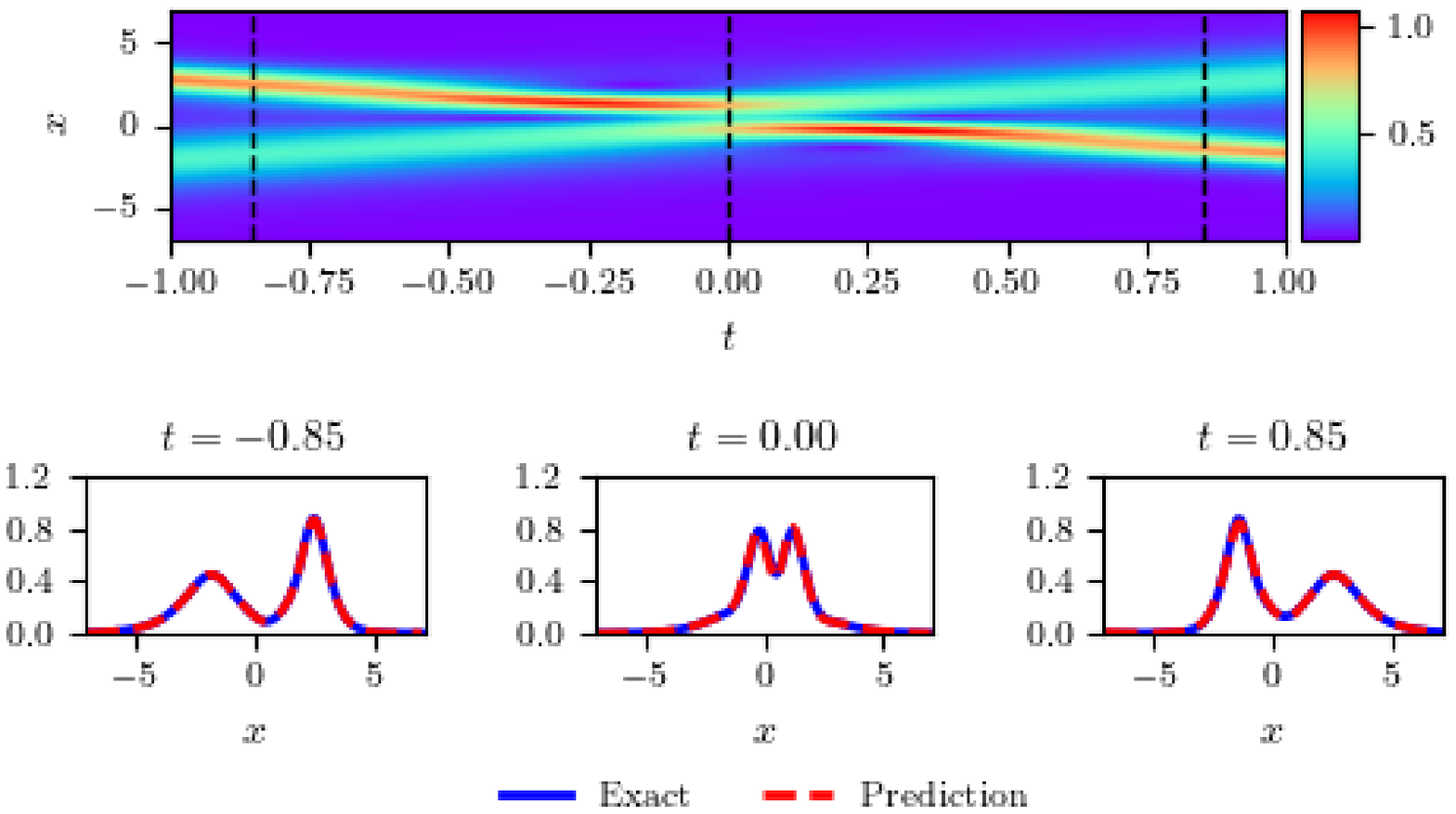}
	{\footnotesize\hspace{3.5cm} (e) RAR-PINN solution $|h_{1}(x,t)|$ \hspace{3.5cm}(f) RAR-PINN solution $|h_{2}(x,t)|$ }  \\
	\caption{The elastic collision two-soliton solutions $h_{1}(x,t)$ and $h_{2}(x,t)$: Comparison among the exact solutions (panels (a) and (b)), predictive  solutions via TPINN (panels (c) and (d)) and predictive results by RAR-PINN (panels (e) and (f)) at different propagation times $ t=-0.85, 0.0 $ and 0.85. The left and right graphs represent the density diagrams and the profiles at different times for the vector two-soliton solutions $h_{1}(x,t)$ and $h_{2}(x,t)$, respectively.
	}
	\label{figure5}
\end{figure}

 To illustrate the collision behavior of vector two-soliton solutions, we take $ [-L,L] $ and $ [-T,T] $ in Eq.~\eqref{equation} as $ [-7,7] $ and $ [-1,1] $, respectively. By choosing $ k_{1}=1+i $, $ k_{2}=2-i $, $ \xi_{1}^{(1)}= \xi_{1}^{(2)}=\xi_{2}^{(1)} = \xi_{2}^{(2)} =1 $, $ \alpha=\beta=2 $ and $ \gamma=0.5+0.5i $ in Solution \eqref{two} the corresponding initial condition is driven as shown in the following formula: $ h_{10}(x)=h_{1}(x,-1) $ and $ h_{20}(x)=h_{2}(x,-1) $. What's more, the boundary condition is expressed as
 \begin{equation}
 	\left\{
 	\begin{aligned}
 		&h_{1}(-7,t)=h_{1}(7,t),h_{2}(-7,t)=h_{2}(7,t),t\in[-1,1], \\
 		&h_{1x}(-7,t)=h_{1x}(7,t),h_{2x}(-7,t)=h_{2x}(7,t).
 	\end{aligned}
 	\right.
 \end{equation}
The computation domain $[-7,7]\times[-1,1]$ is divided into $[300\times201]$ data points, and then vector two-soliton solutions $h_{1}$ and $h_{2}$ are discretized into 201 snapshots on regular space-time grid with $\bigtriangleup t=0.01$. A training data set containing initial data and boundary data are generated by randomly sampling \cite{1994Random}. Specifically, the number of boundary points is $N_{b}=100$, and the number of initial points is $N_{0}=100$. Firstly, we choose 4000 points selected as the residual points and then add 15 more new points adaptively via RAR with $m=3$ and $\varepsilon_{0}=0.055$. In order to compare these two algorithms, we select $N_{0}=100$, $N_{b}=100$ and $N_{f}=4015$ to obtain the vector two-soliton solutions using the TPINN algorithm. The predicted vector two-soliton solutions $h_{1}(x,t)$, and $h_{2}(x,t)$ can be obtained by minimizing the loss function. Here a fully connected neural network with 6 hidden layers and 32 neurons per hidden layer is used to learn the parameters of the neural network. Besides, Figure \ref{figure5} specifically demonstrates the effectiveness of RAR-PINN in solving the elastic collision vector two-soliton solutions.

As shown in Figure \ref{figure5}, we can clearly see the reconstructed vector two-soliton space-time dynamics. The density plots with the corresponding peak scale for diverse dynamics include exact dynamics and two dynamics learned by different methods. More specifically, Figure \ref{figure5} (a-b) are the exact interactions of two solitons, while Figure \ref{figure5} (c-d) and Figure \ref{figure5} (e-f) show the results given by TPINN and RAR-PINN methods, respectively. The blue and red curves represent the exact and predicted solutions, respectively. The predicted solutions obtained by the two methods are compared with the exact solutions, respectively. The relative $ \mathbb{L}_{2} $ errors of the TPINN model are $ 7.381569\times10^{-2} $ for $ h_{1}(x,t) $ and $ 7.395515\times10^{-2} $ for $ h_{2}(x,t) $; while the relative $ \mathbb{L}_{2} $ errors of the RAR-PINN model are $ 3.764622\times10^{-2} $ for $ h_{1}(x,t) $ and $ 3.756750\times10^{-2} $ for $ h_{2}(x,t) $.


\begin{figure}[htbp]
	\centering
	\includegraphics[scale=0.06]{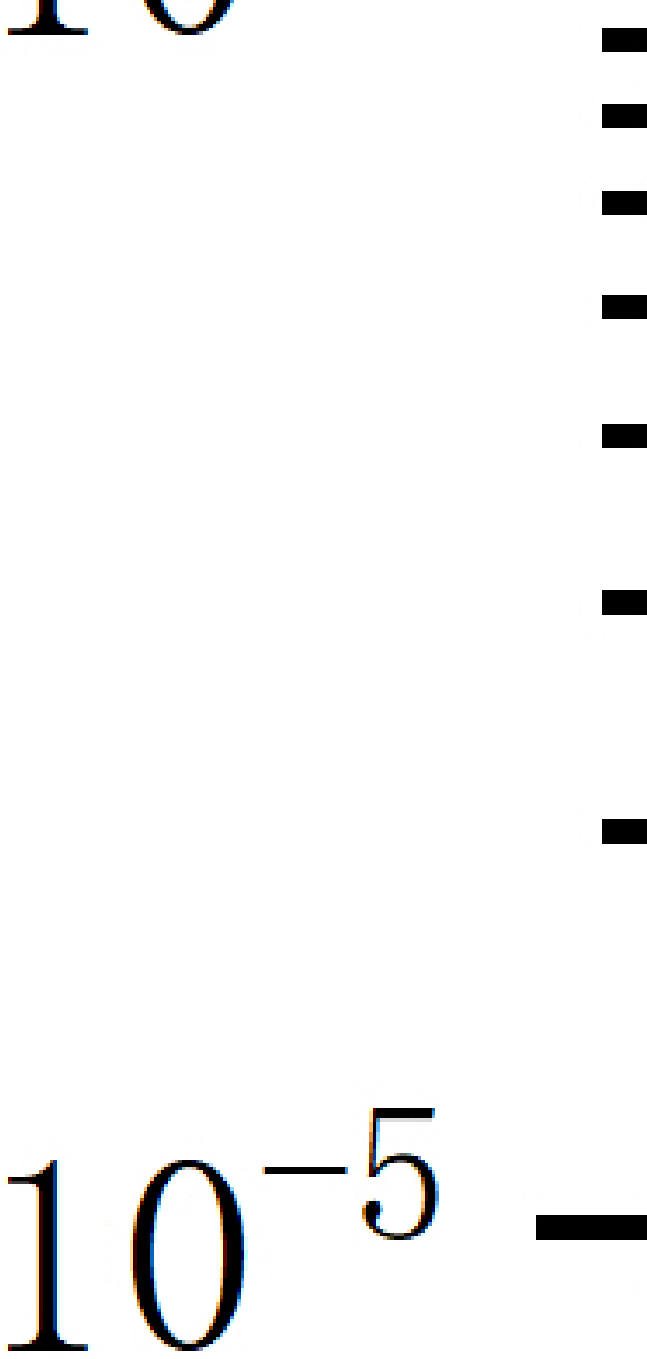}\hspace{1cm}
	\includegraphics[scale=0.06]{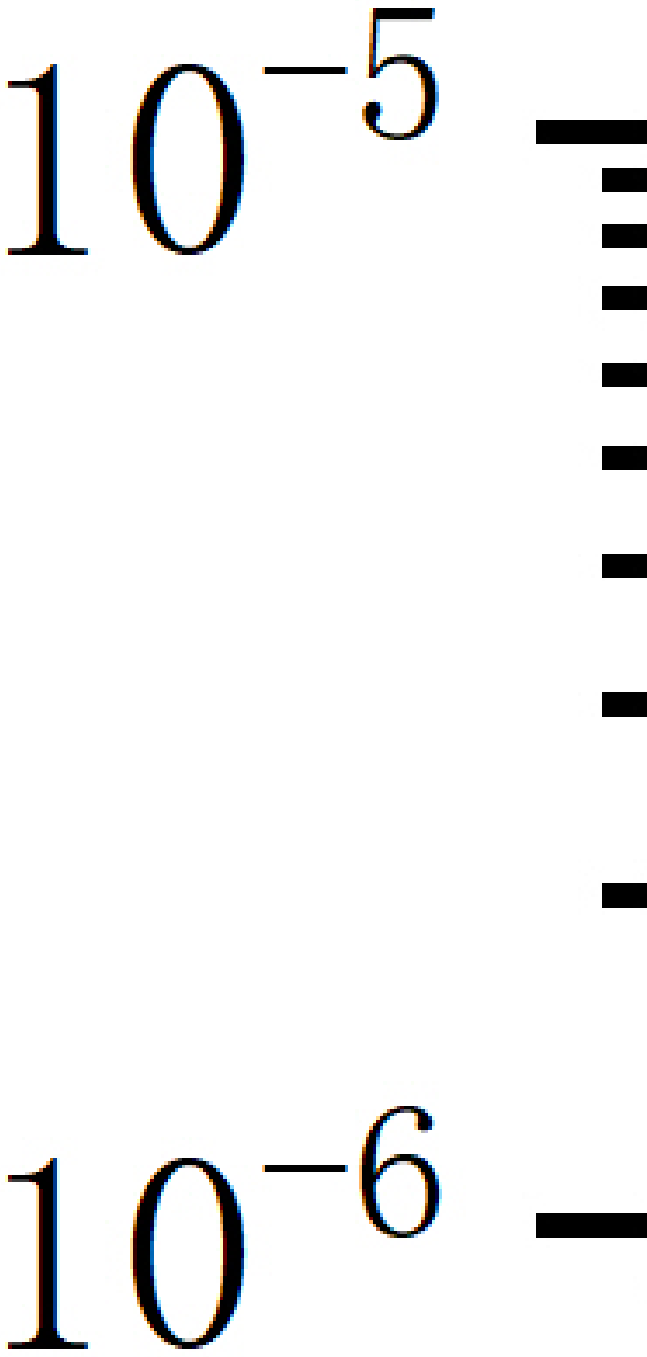}
	{\footnotesize\hspace{9cm} (a) \hspace{8cm}(b) } \\
	\includegraphics[scale=0.06]{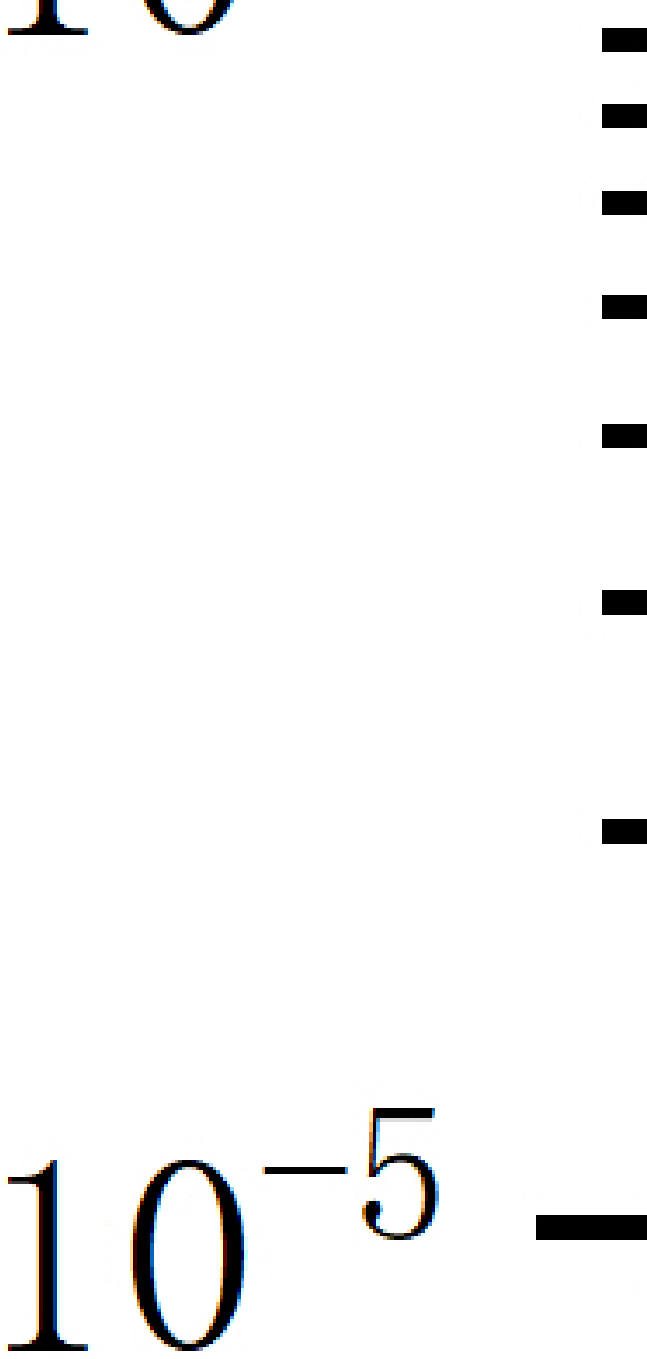}\hspace{1cm}
	\includegraphics[scale=0.06]{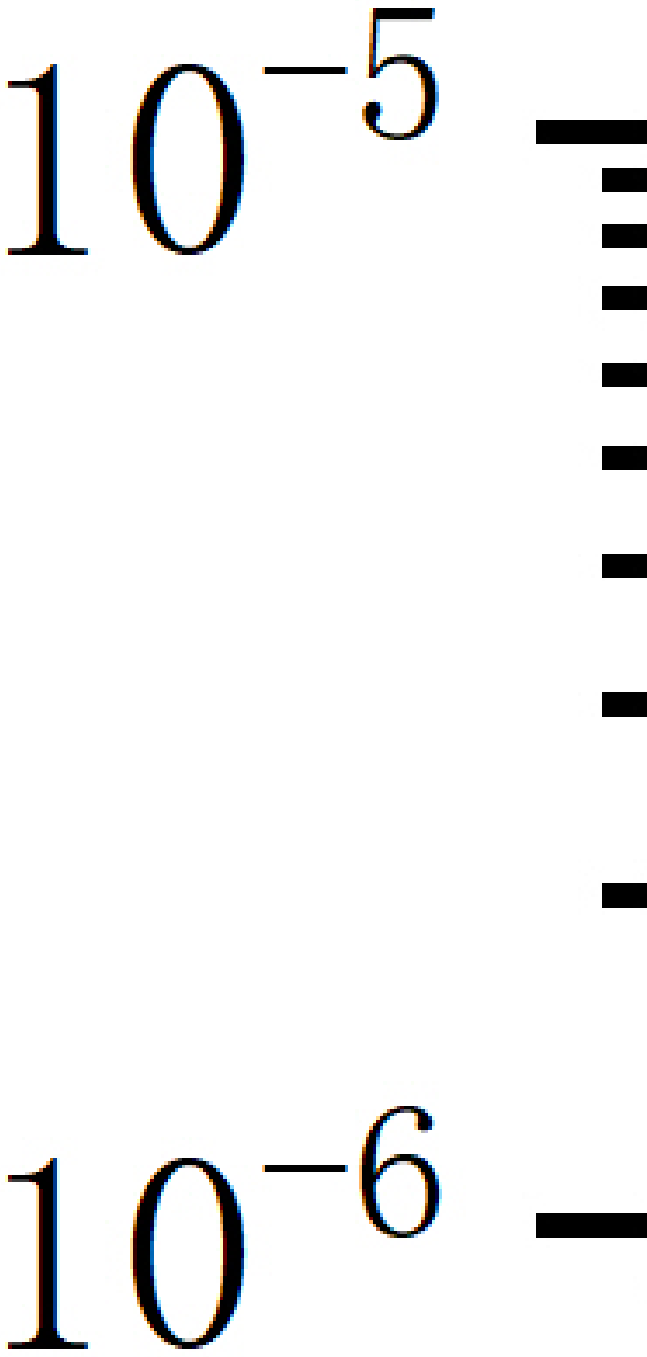}
	{\footnotesize\hspace{9cm} (c) \hspace{8cm}(d) } \\  	
	\caption{The loss function curve figures of vector two-soliton solutions $h_{1}(x,t)$ and $h_{2}(x,t)$ arising from (a) TPINN algorithm with the 10000 iterations Adam optimization; (b) TPINN algorithm with the 9423 iterations L-BFGS optimization; (c) RAR-PINN algorithm with the 10000 iterations Adam optimization; (d) RAR-PINN algorithm with the 13272 iterations L-BFGS optimization.}
	\label{figure6}
\end{figure}

During the process of neural network training, we first use Adam optimizer \cite{diederik2014adam} to train 10,000 times, then use L-BFGS optimizer \cite{1989On} to train. Adam is an algorithm for optimizing stochastic objective functions based on first-order gradients. The main advantages of the Adam method are: simple and easy adaptive estimation based on low order moments on the gradient; computationally efficient and low memory requirements well suited for problems with large data and/or parameters and suitable for non-smooth objectives and problems with very noisy and/or sparse gradients \cite{diederik2014adam}. Figure \ref{figure6}(a) shows the loss function curve using Adam's algorithm trained 10,000 iterations and Figure \ref{figure6}(b) represents the loss function descent curve using L-BFGS's algorithm trained 9423 iterations. More specially, initial condition error is $ loss_{0}=1.401873\times10^{-5}$; boundary condition error is $ loss_{b}=2.5388058\times10^{-6}$; the collocation points error is $ loss_{f}=5.429575\times10^{-5} $ and overall error is $ Loss= 7.085329\times10^{-5}$.
Figure \ref{figure6}(c) and (d) represent the loss function curves with RAR-PINN using Adam optimizer and L-BFGS optimizer, respectively. After training, initial condition error is $ loss_{0}=6.1637432\times10^{-6}$; boundary condition error is $ loss_{b}=1.3867134\times10^{-6}$; the collocation points error is $ loss_{f}=2.861435\times10^{-5} $ and overall error is $ Loss= 3.616529\times10^{-5}$.

\begin{figure}[htbp]
	\centering
	\includegraphics[scale=0.06]{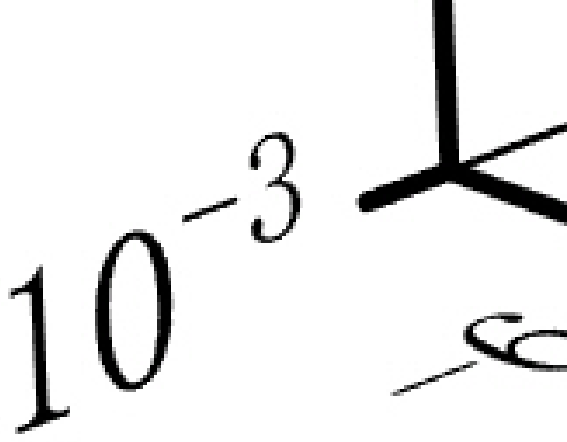}\hspace{1cm}
	\includegraphics[scale=0.06]{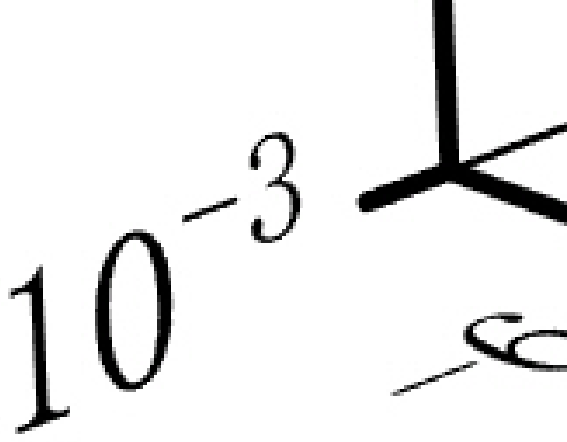}
	{\footnotesize\hspace{9cm} (a) \hspace{8cm}(b) } 	
	\caption{Three-dimensional stereogram of residual error: (a) three-dimensional stereogram residual error without adding adaptive points; (b) three-dimensional stereogram residual error after adding adaptive points.
	}
	\label{figure7}
\end{figure}

 The three-dimensional diagram of the residual error is shown in Figure \ref{figure7}. Figure \ref{figure7}(a) is three-dimensional stereogram residual error without adding adaptive points, from which we can see that the residuals of some positions before the adaptive pointing are relatively large up to 2.270. Figure \ref{figure7}(b) denotes the denotes the three-dimensional stereogram residual error after adding adaptive points, which indicates that the residuals of all positions after the adaptive pointing can be controlled within 1.30. From these experimental results, we can find that the RAR-PINN method is more efficient than TPINN for solving the elastic collision of vector two-soliton solutions.

\subsubsection{Inelastic collision of vector two-soliton solution}

In this section, by choosing the parameters  $k_{1}=1+i$, $k_{2}=2-i$, $\xi_{1}^{(1)}= \xi_{1}^{(2)}=\xi_{2}^{(2)} =1$, $\xi_{2}^{(1)}=\dfrac{39+80i}{89}$, $\alpha=\beta=2$ and $\gamma=0.5+0.5i$, we get the soliton collision behaviors. We take $[-L,L]$ and $[-T,T]$ in System (1) as $[-4,4]$ and $[-0.3,0.3]$, respectively. The initial values can be computed by Eq.~(\ref{two}) at $t=-0.3$ and the boundary values are periodic boundary condition. We divide the domain into $[400\times301]$ data points. The vector two-soliton solutions $h_{1}$ and $h_{2}$ are discretized into 301 snapshots. We predict the data-driven solution based on the neural networks architecture with 6 layers and 32 neurons per hidden layer. The training set of the TPINN model consists of initial training points $ N_{0} = 150 $, boundary training points $ N_{b}=150 $ and collocation points $ N_{f}=15040 $. As for the RAR-PINN model, we first choose $ N_{0} = 150 $, $ N_{b}=150 $ and $ N_{f}=15000 $, then add 100 more collocation points adaptivly. In particular, we choose the threshold $ \varepsilon_{0}=0.13 $ and the number of points added each time $m=10$.


\begin{figure}[htbp]
	\centering
	\includegraphics[scale=0.4]{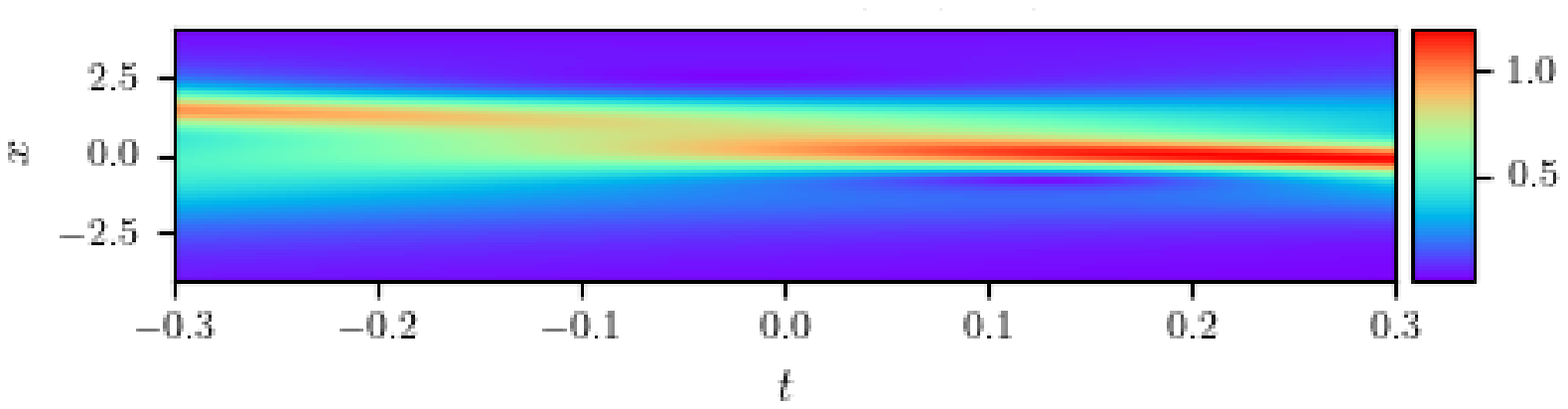}\hspace{1cm}
	\includegraphics[scale=0.4]{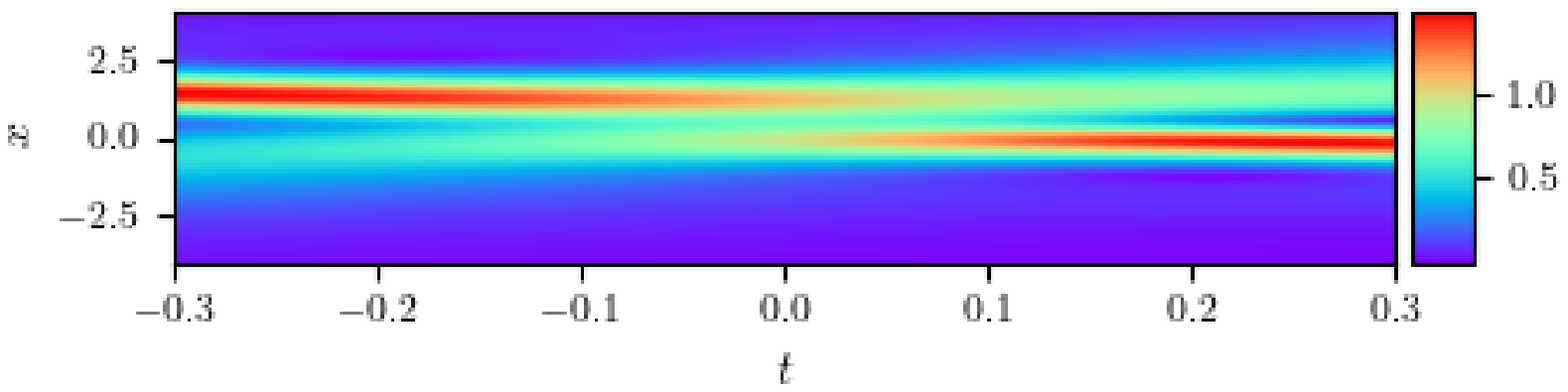}
	{\footnotesize\hspace{4cm} (a) Exact solution $ |h_{1}(x,t)| $ \hspace{4.5cm}(b) Exact solution $ |h_{2}(x,t)| $ } \\
	\vspace{5mm}
	\includegraphics[scale=0.4]{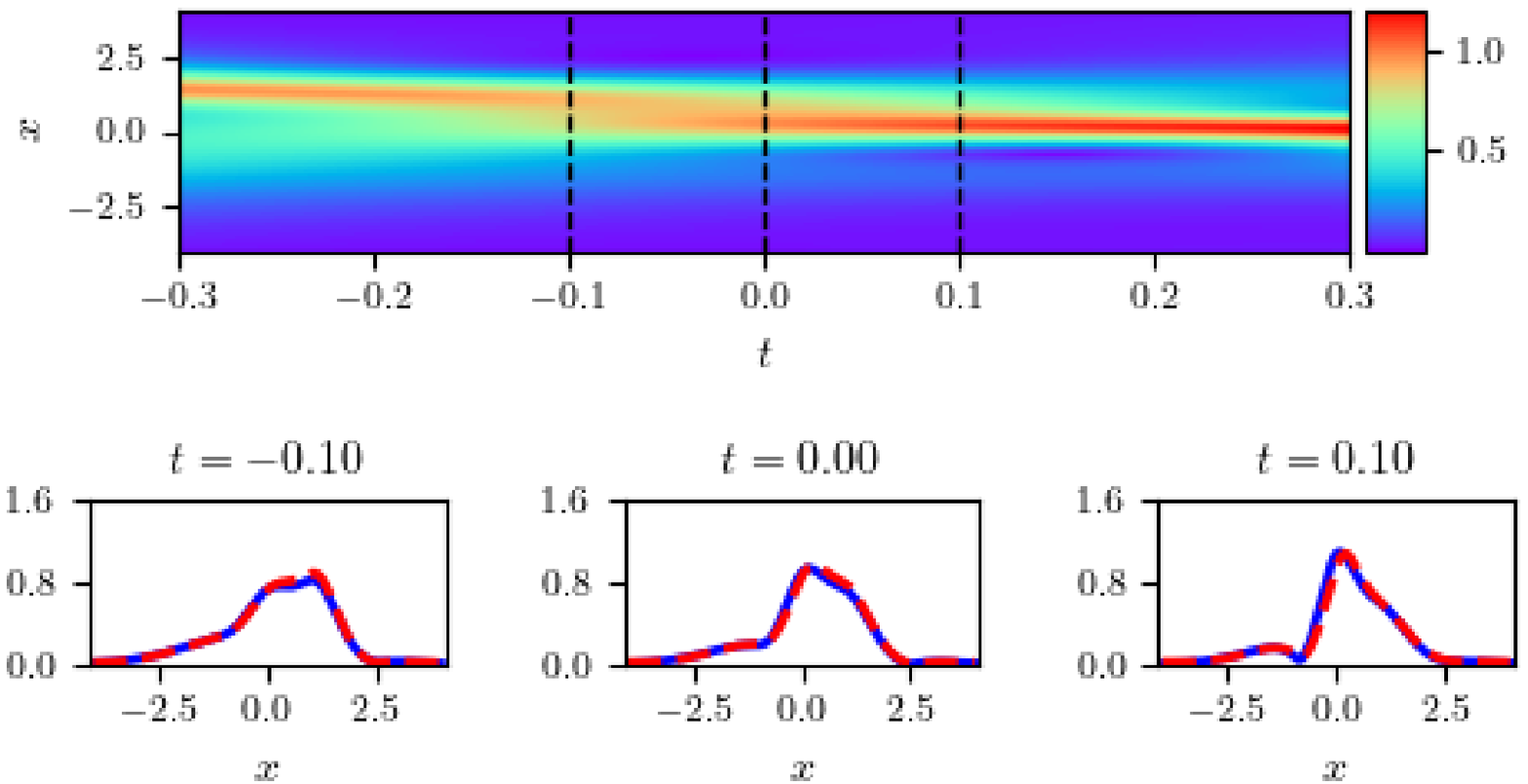}\hspace{1cm}
	\includegraphics[scale=0.4]{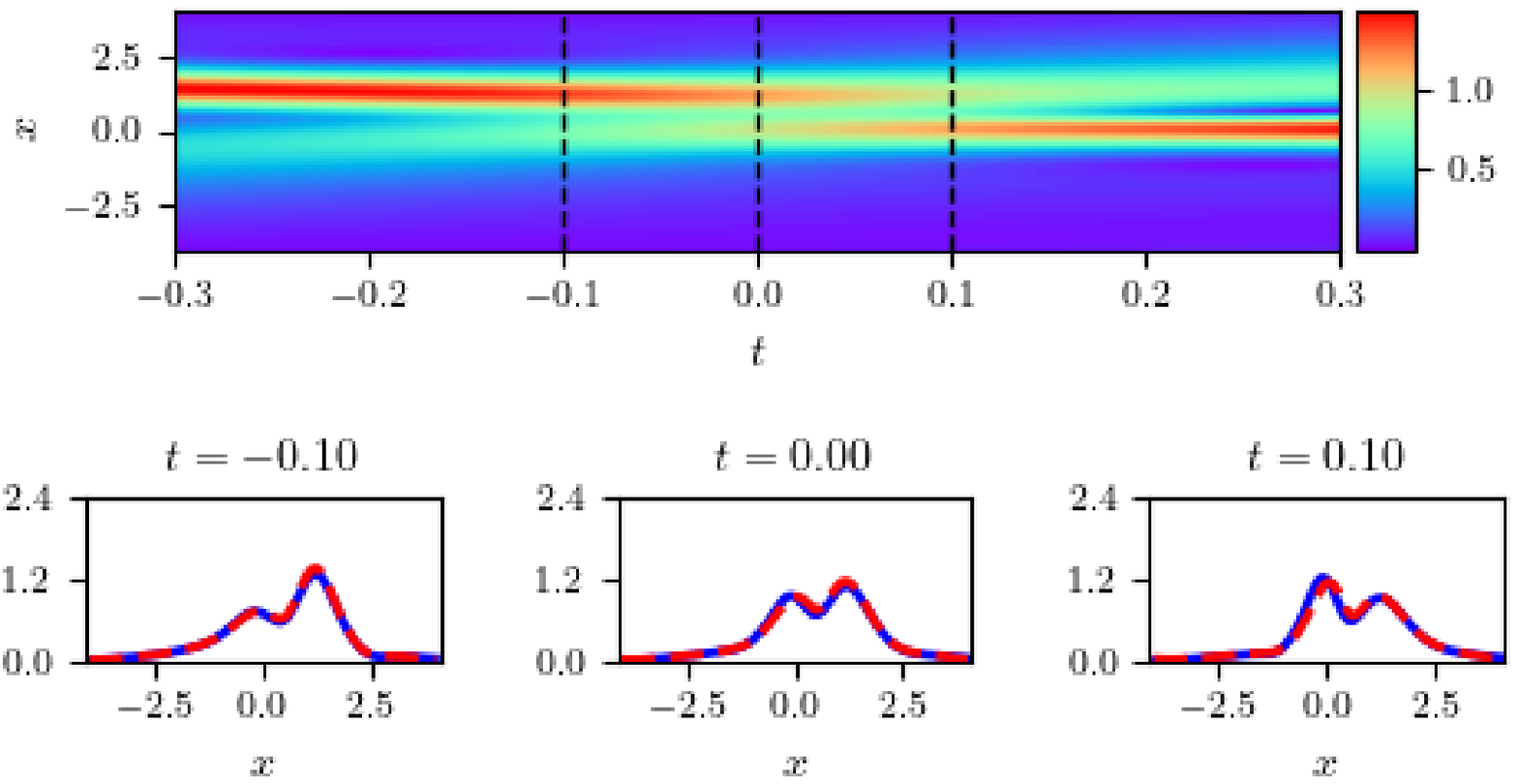}
	{\footnotesize\hspace{3cm} (c) Tractional PINN solution $ |h_{1}(x,t)| $ \hspace{3cm}(d) Tractional PINN solution $ |h_{2}(x,t)| $ } \\
	\vspace{5mm}
	\includegraphics[scale=0.4]{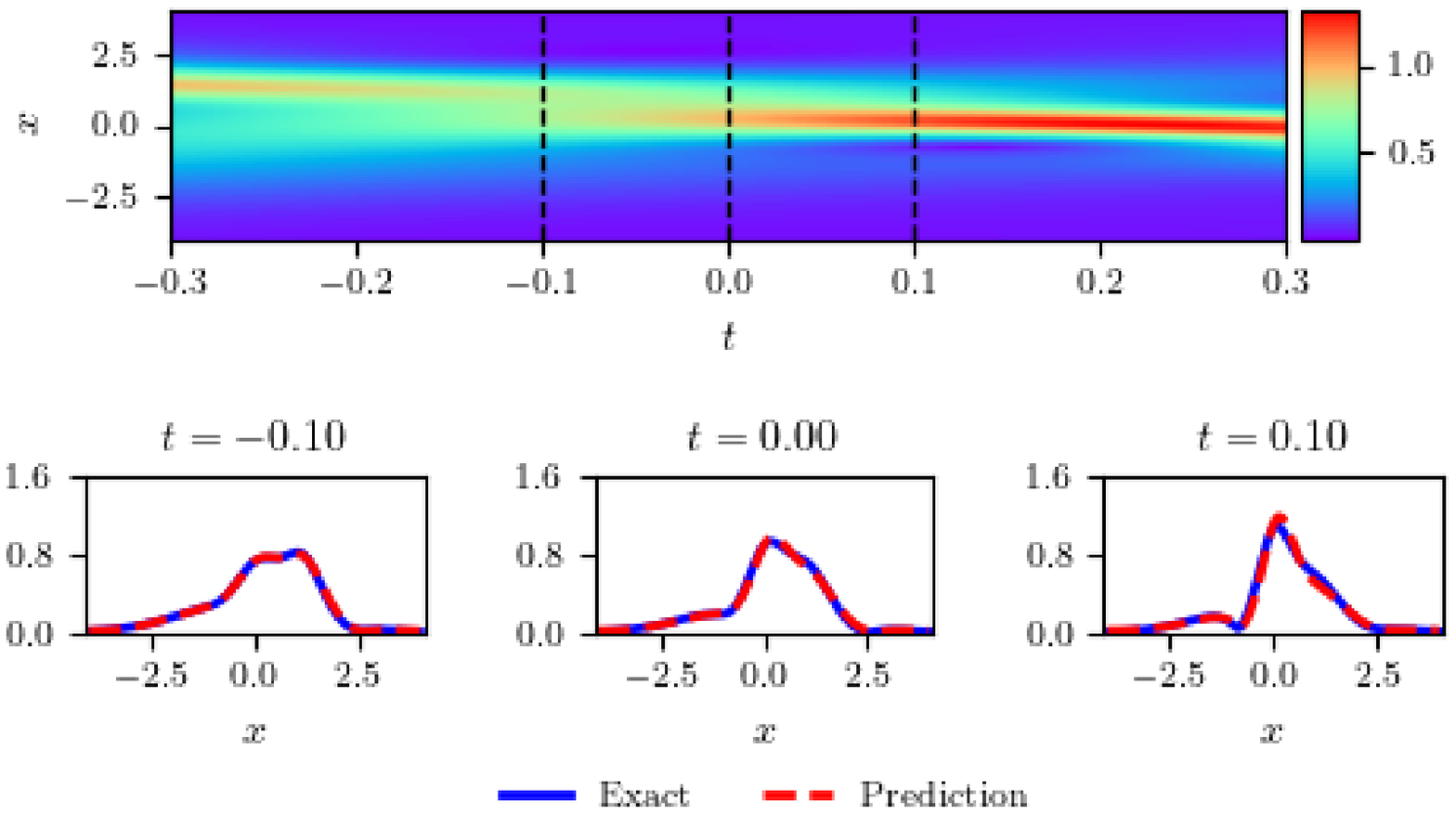}\hspace{1cm}
	\includegraphics[scale=0.4]{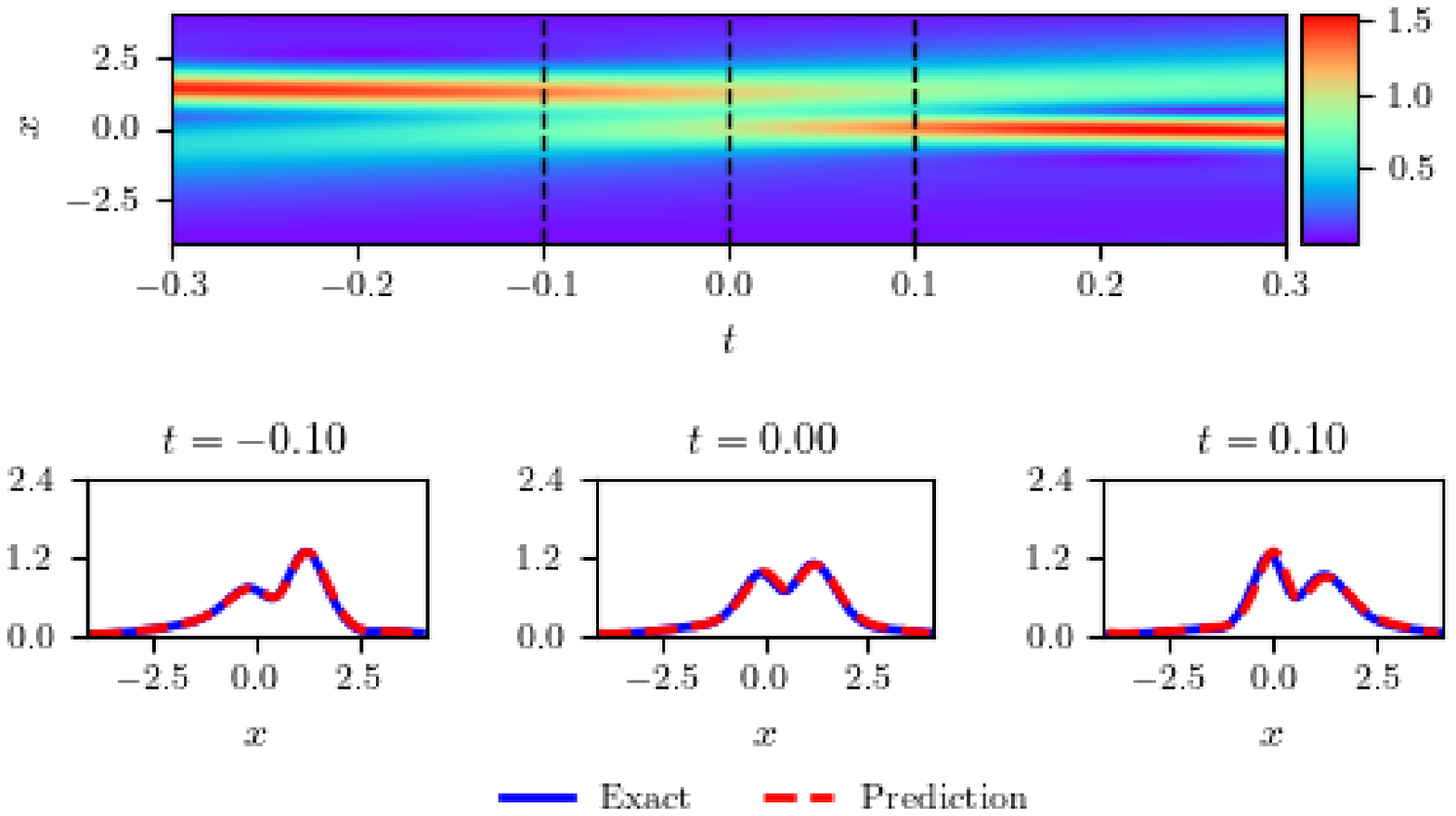}
	{\footnotesize\hspace{3.5cm} (e) RAR-PINN solution $ |h_{1}(x,t)| $ \hspace{3.5cm}(f) RAR-PINN solution $ |h_{2}(x,t)| $ }  \\
	\caption{The shape-changing vector two-soliton solutions $ h_{1}(x,t) $ and $ h_{2}(x,t) $: Comparison among the exact solutions (panels (a) and (b)), predictive solutions via TPINN (panels (c) and (d)) and predictive results by RAR-PINN (panels (e) and (f)) at different propagation times $ t=-0.1, 0.0 $ and 0.1. The left and right graphs represent the density diagrams and the profiles at different times for the shape-changing vector two-soliton solutions $ h_{1}(x,t) $ and $ h_{2}(x,t) $, respectively.
	}
	\label{figure8}
\end{figure}

\begin{figure}[htbp]
	\centering
	\includegraphics[scale=0.06]{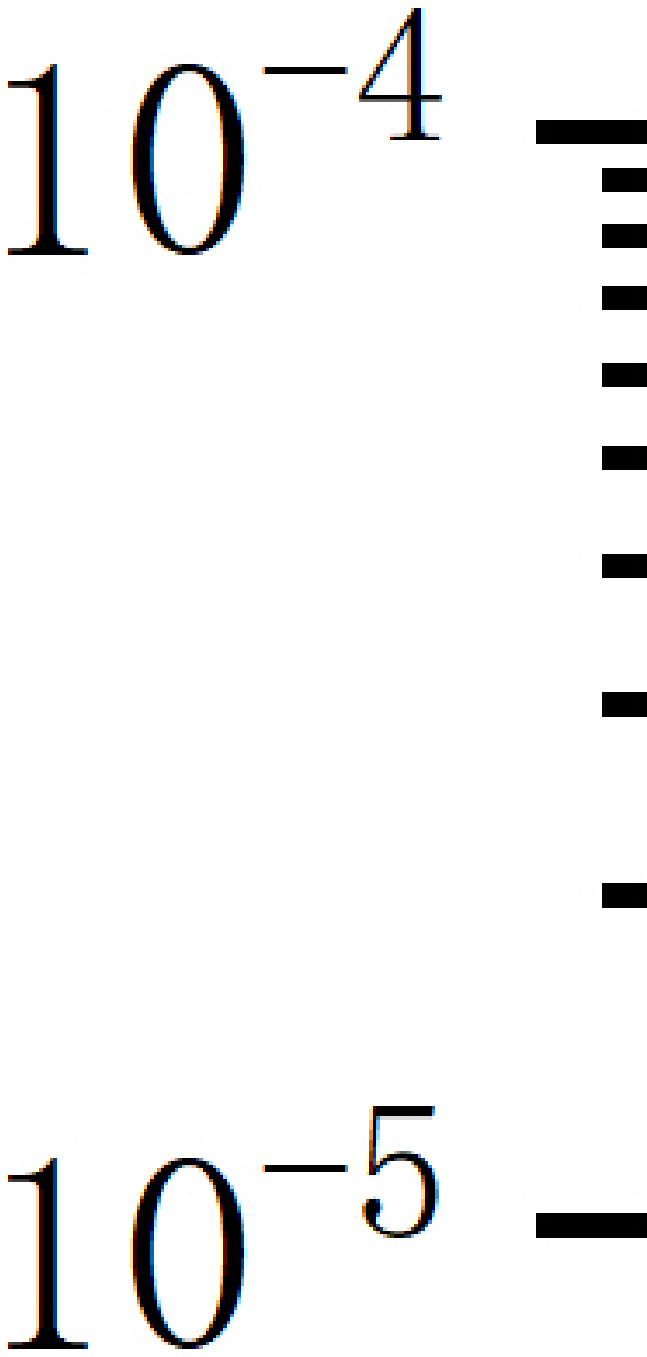}\hspace{1cm}
	\includegraphics[scale=0.06]{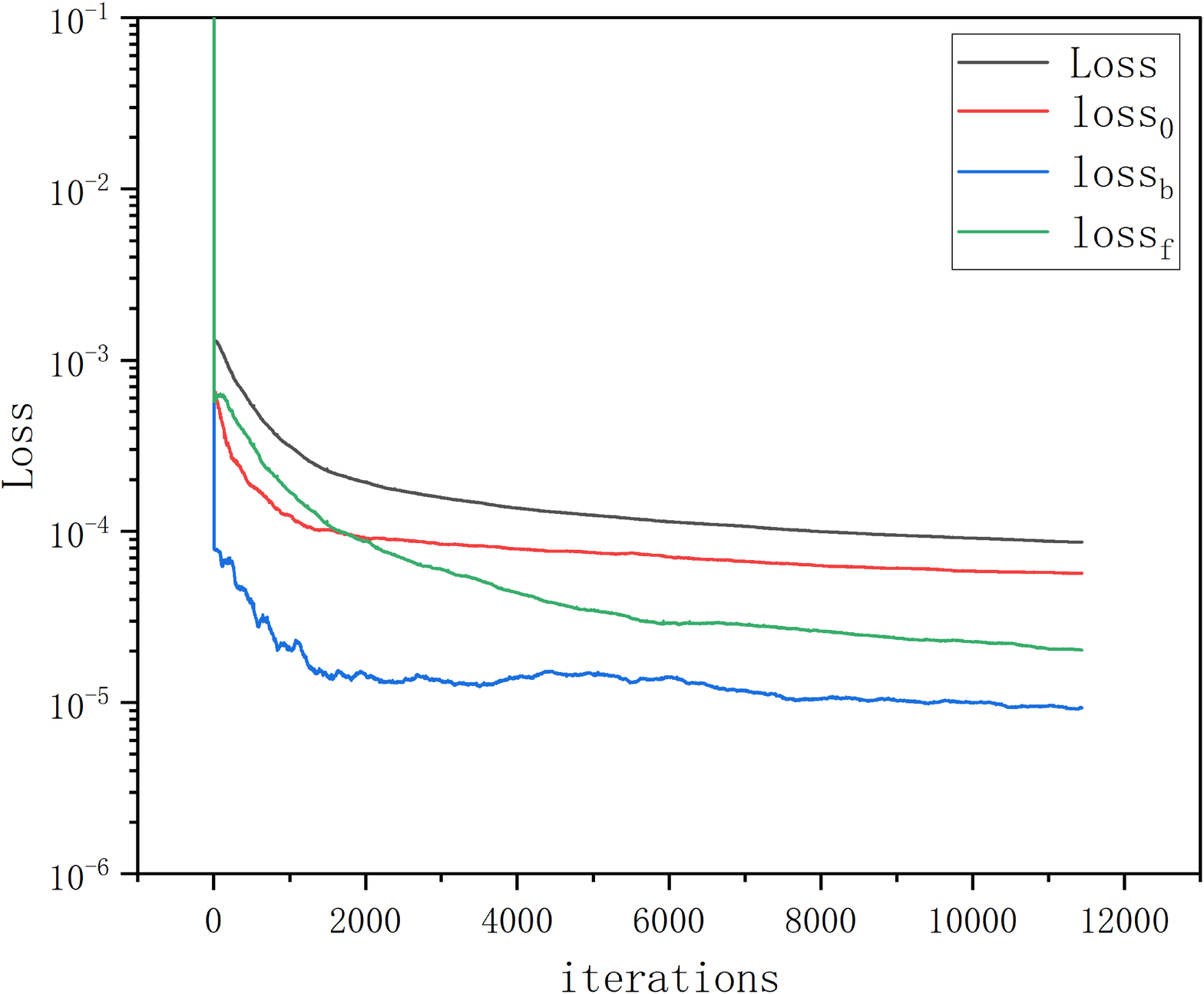}
	{\footnotesize\hspace{9cm} (a) \hspace{8cm}(b) } \\
	\includegraphics[scale=0.06]{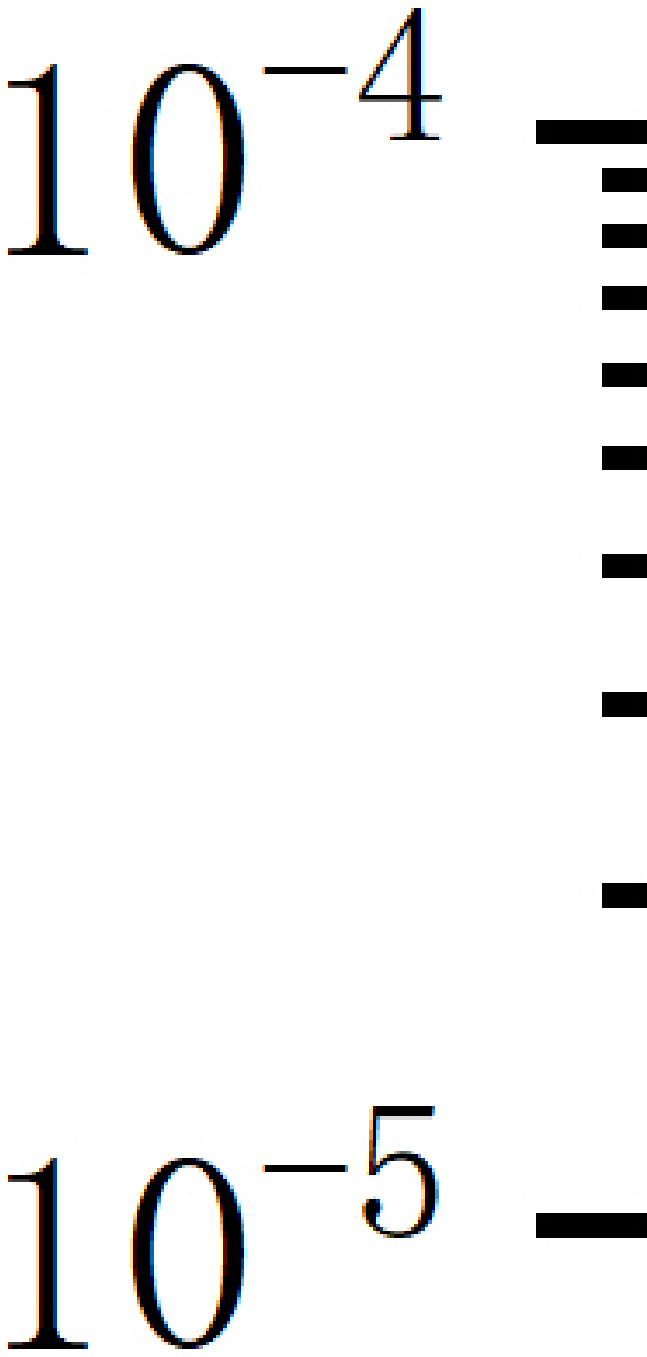}\hspace{1cm}
	\includegraphics[scale=0.06]{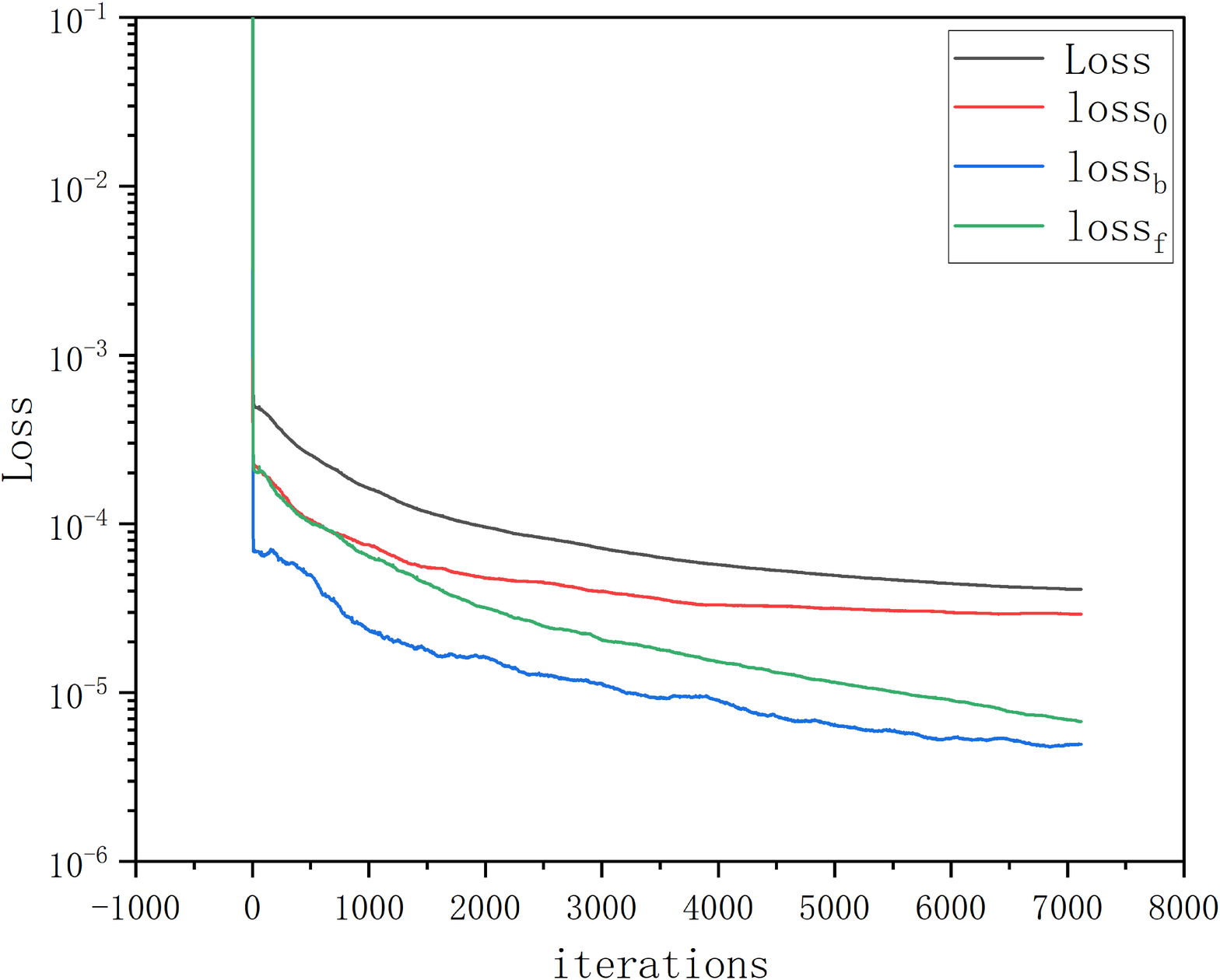}
	{\footnotesize\hspace{9cm} (c) \hspace{8cm}(d) } \\  	
	\caption{The loss function curve figures of shape-changing vector two-soliton solutions $h_{1}(x,t)$ and $h_{2}(x,t)$ arising from (a) TPINN algorithm with the 20000 iterations Adam optimization; (b) TPINN algorithm with the 11438 iterations L-BFGS optimization; similarly, (c)RAR-PINN algorithm with the 20000 iterations Adam optimization; (d) RAR-PINN algorithm with the 7117 iterations L-BFGS optimization.}
	\label{figure9}
\end{figure}

\begin{figure}[htbp]
	\centering
	\includegraphics[scale=0.06]{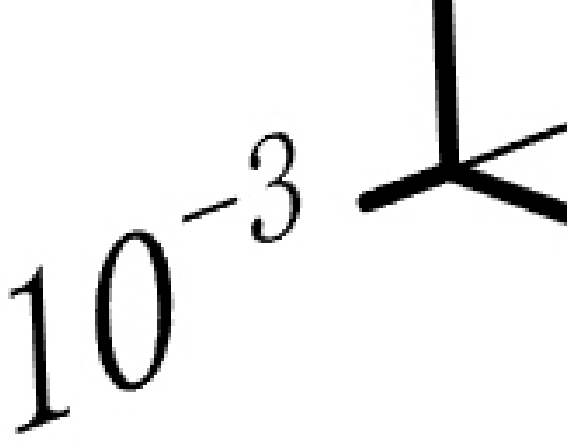}\hspace{1cm}
	\includegraphics[scale=0.06]{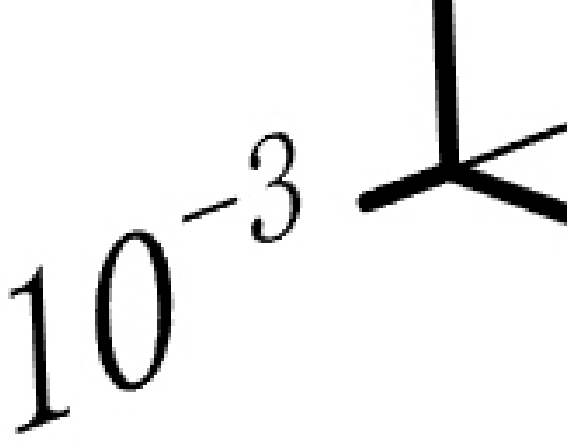}
	{\footnotesize\hspace{9cm} (a) \hspace{8cm}(b) } 	
	\caption{Three-dimensional stereogram of residual error: (a) three-dimensional stereogram residual error without adding adaptive points; (b) three-dimensional stereogram residual error after adding adaptive points.
	}
	\label{figure10}
\end{figure}

 Figure \ref{figure8} demonstrates the inelastic collision of vector two-soliton solution predicted by the TPINN model and the RAR-PINN model, respectively. The relative $\mathbb{L}_{2}$ errors of the shape-changing vector two-soliton solutions $ h_{1}(x,t) $ and $ h_{2}(x,t) $ via the TPINN model are $1.196291\times10^{-1}$ and $1.205316\times10^{-1}$, respectively; while the relative $\mathbb{L}_{2}$ errors via the RAR-PINN model are $8.64846\times10^{-2}$ for $h_{1}(x,t)$ and $8.71967\times10^{-2}$ for $h_{2}(x,t)$. Adam optimizer and L-BFGS optimizer are used to minimize the loss function. Figure \ref{figure9}(a) demonstrates the loss function descent curve of the TPINN model with the Adam optimizer for 20,000 iterations. Figure \ref{figure9}(b) represents the loss function curve of the TPINN model with the L-BFGS optimizer for 11438 iterations. We can find the initial error $loss_{0}=5.6866372 \times10^{-5}$, the boundary error $loss_{b}=9.3172368\times10^{-6}$, the residual error $loss_{f}=2.0258679\times10^{-5}$ and the total error $Loss=8.644289\times10^{-5}$. Figure \ref{figure9}(c) indicates the loss function curve of the RAR-PINN model via the Adam optimizer for 20,000 iterations. Figure \ref{figure9}(d) illustrates the loss function curve of the RAR-PINN model with the L-BFGS optimizer for 7117 iterations. Furthermore, we can observe the initial error $ loss_{0}=2.923557\times10^{-5}$, the boundary error $ loss_{b}=4.950383\times10^{-6} $, the residual $ loss_{f}=6.765904\times10^{-6} $ and the total error $ Loss=4.095181\times10^{-5} $. On the whole, the RAR-PINN model is more efficient than the TPINN model in both relative and absolute errors.
 Figure \ref{figure10}(a) and (b) demonstrate the residual errors of the three-dimensional stereogram without the addition of adaptive points and the residual errors of the three-dimensional stereogram after the addition of adaptive points, respectively.
 Through the above comparison, we can see that the RAR-PINN model is more efficient for the shape-changing soliton collisions of Eq.~\eqref{equation}.

\subsection{Vector three-soliton solutions}

The exact vector three-soliton solution of the CGNLS equation has been studied by various methods, that is given as\cite{VISHNUPRIYA2016366}
\begin{equation}
	\begin{aligned}
		h_{1}=&\frac{\xi_{1}^{(1)}e^{\eta_{1}}+\xi_{2}^{(1)}e^{\eta_{2}}+\xi_{3}^{(1)}e^{\eta_{3}}+e^{\eta_{1}+\eta_{1}^{*}+\eta_{2}+\delta_{11}}+e^{\eta_{1}+\eta_{1}^{*}+\eta_{3}+\delta_{21}}+e^{\eta_{1}+\eta_{2}^{*}+\eta_{2}+\delta_{21}}}{D}\\
		&+\frac{e^{\eta_{2}+\eta_{2}^{*}+\eta_{3}+\delta_{41}}+e^{\eta_{3}+\eta_{3}^{*}+\eta_{1}+\delta_{51}}+e^{\eta_{3}+\eta_{3}^{*}+\eta_{2}+\delta_{61}}+e^{\eta_{1}^{*}+\eta_{2}+\eta_{3}+\delta_{71}+}e^{\eta_{1}+\eta_{2}^{*}+\eta_{3}+\delta_{81}}}{D}\\
		&+\frac{e^{\eta_{1}+\eta_{2}+\eta_{3}^{*}+\delta_{91}}+e^{\eta_{1}+\eta_{1}^{*}+\eta_{2}+\eta_{2}^{*}+\eta_{3}+\tau_{11}}+e^{\eta_{1}+\eta_{1}^{*}+\eta_{3}+\eta_{3}^{*}+\eta_{2}+\tau_{21}}+e^{\eta_{2}+\eta_{2}^{*}+\eta_{3}+\eta_{3}^{*}+\eta_{1}+\tau_{31}}}{D},\\
		h_{2}=&\frac{\xi_{1}^{(2)}e^{\eta_{1}}+\xi_{2}^{(2)}e^{\eta_{2}}+\xi_{3}^{(2)}e^{\eta_{3}}+e^{\eta_{1}+\eta_{1}^{*}+\eta_{2}+\delta_{12}}+e^{\eta_{1}+\eta_{1}^{*}+\eta_{3}+\delta_{22}}+e^{\eta_{1}+\eta_{2}^{*}+\eta_{2}+\delta_{32}}}{D}\\
		&+\frac{e^{\eta_{2}+\eta_{2}^{*}+\eta_{3}+\delta_{42}}+e^{\eta_{3}+\eta_{3}^{*}+\eta_{1}+\delta_{52}}+e^{\eta_{3}+\eta_{3}^{*}+\eta_{2}+\delta_{62}}+e^{\eta_{1}^{*}+\eta_{2}+\eta_{3}+\delta_{72}+}e^{\eta_{1}+\eta_{2}^{*}+\eta_{3}+\delta_{82}}}{D}\\
		&+\frac{e^{\eta_{1}+\eta_{2}+\eta_{3}^{*}+\delta_{92}}+e^{\eta_{1}+\eta_{1}^{*}+\eta_{2}+\eta_{2}^{*}+\eta_{3}+\tau_{12}}+e^{\eta_{1}+\eta_{1}^{*}+\eta_{3}+\eta_{3}^{*}+\eta_{2}+\tau_{22}}+e^{\eta_{2}+\eta_{2}^{*}+\eta_{3}+\eta_{3}^{*}+\eta_{1}+\tau_{32}}}{D},\label{intwo}
	\end{aligned}
\end{equation}
where
\begin{equation}
	\begin{split}
		\begin{aligned}
			D=&1+e^{\eta_{1}+\eta_{1}^{*}+R_{1}}+e^{\eta_{2}+\eta_{2}^{*}+R_{2}}+e^{\eta_{3}+\eta_{3}^{*}+R_{3}}+e^{\eta_{1}+\eta_{2}^{*}+\delta_{10}}+e^{\eta_{1}^{*}+\eta_{2}+\delta_{10}^{*}}\\
			&+e^{\eta_{1}+\eta_{3}^{*}+\delta_{20}}+e^{\eta_{1}^{*}+\eta_{3}+\delta_{20}^{*}}+e^{\eta_{2}+\eta_{3}^{*}+\delta_{30}}+e^{\eta_{2}^{*}+\eta_{3}+\delta_{30}^{*}}+e^{\eta_{1}+\eta_{1}^{*}+\eta_{2}+\eta_{2}^{*}+R_{4}}\\
			&+e^{\eta_{1}+\eta_{1}^{*}+\eta_{3}+\eta_{3}^{*}+R_{5}}+e^{\eta_{2}+\eta_{2}^{*}+\eta_{3}+\eta_{3}^{*}+R_{6}}+e^{\eta_{1}+\eta_{1}^{*}+\eta_{2}+\eta_{3}^{*}+\tau_{10}}+e^{\eta_{1}+\eta_{1}^{*}+\eta_{3}+\eta_{2}^{*}+\tau_{10}^{*}}\\
			&+e^{\eta_{1}+\eta_{1}^{*}+\eta_{1}+\eta_{1}^{*}+\eta_{3}+\eta_{3}^{*}+R_{7}},\\
			\eta_{l}=&k_{l}(x+ik_{l}t), l=1,2,3.		
		\end{aligned}
	\end{split}
\end{equation}
The expressions for the other parameters of the vector three-soliton solutions have been given in detail in the Ref.\cite{VISHNUPRIYA2016366}.
We make choose the vector three-soliton solutions $h_{1}(x,t)$ and $h_{2}(x,t)$ with the parameter choices $k_{1}=1+i$, $ k_{2}=2-i$, $\xi_{1}^{(1)}=\xi_{1}^{(2)}=\xi_{2}^{(1)}=\xi_{2}^{(2)}=1$, $\alpha=\beta=2$ and $\gamma=0.5+0.5i$. Besides, we take $[-L,L]$ and $[-T,T]$ in Eq.~\eqref{equation} as $[-6,6]$ and $[-0.8,0.8]$, respectively. The corresponding initial condition is driven as $h_{10}=h_{1}(x,-1)$ and $h_{20}=h_{2}(x,-1)$. Similarly, here we still use periodic boundary condition to numerically simulate the vector three-soliton solutions of Eq.~\eqref{equation}.The computation domain $[-6,6]\times[-0.8,0.8]$ is divided into $[400\times301]$ data points, three-soliton solutions $h_{1}$ and $h_{2}$ are discretized into 301 snapshots on regular space-time grid with $ \bigtriangleup t=0.005 $. A training data set containing initial data and boundary data are generated by randomly sampling \cite{1994Random}. Specifically, the number of boundary points is $N_{b}=100$ and the number of initial points is $N_{0}=120$. Firstly, we choose 6000 points selected as the residual points and then add 18 more new points adaptively via RAR with $m=3$ and $\varepsilon_{0}=0.07$. In order to compare these two algorithms, we select $N_{0}=120$, $N_{b}=100$ and $N_{f}=6018$ to obtain the vector three-soliton solutions using the TPINN algorithm. The predicted vector three-soliton solutions $h_{1}(x,t)$ and $h_{2}(x,t)$ can be obtained by minimizing the loss function.
\begin{figure}[htbp]
	\centering
	\includegraphics[scale=0.4]{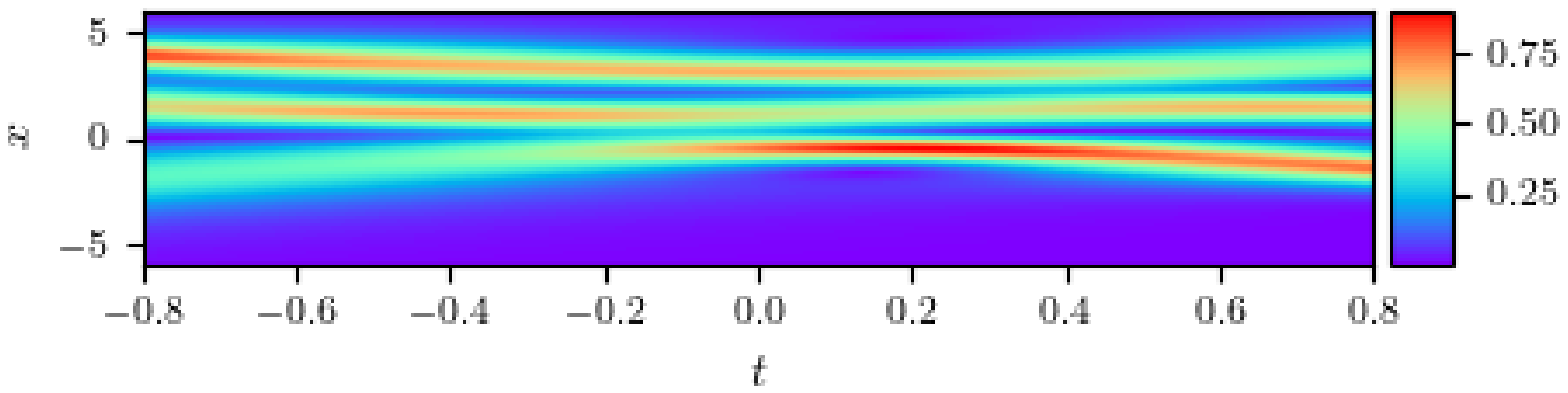}\hspace{1cm}
	\includegraphics[scale=0.4]{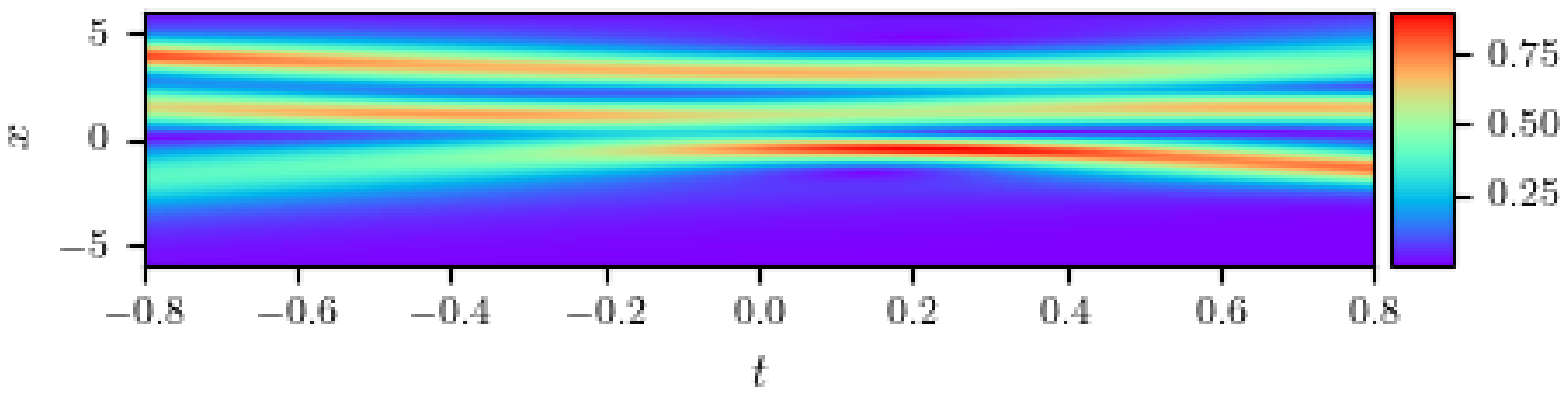}
	{\footnotesize\hspace{4cm} (a) Exact solution $ |h_{1}(x,t)| $ \hspace{4.5cm}(b) Exact solution $ |h_{2}(x,t)| $ } \\
	\vspace{5mm}
	\includegraphics[scale=0.4]{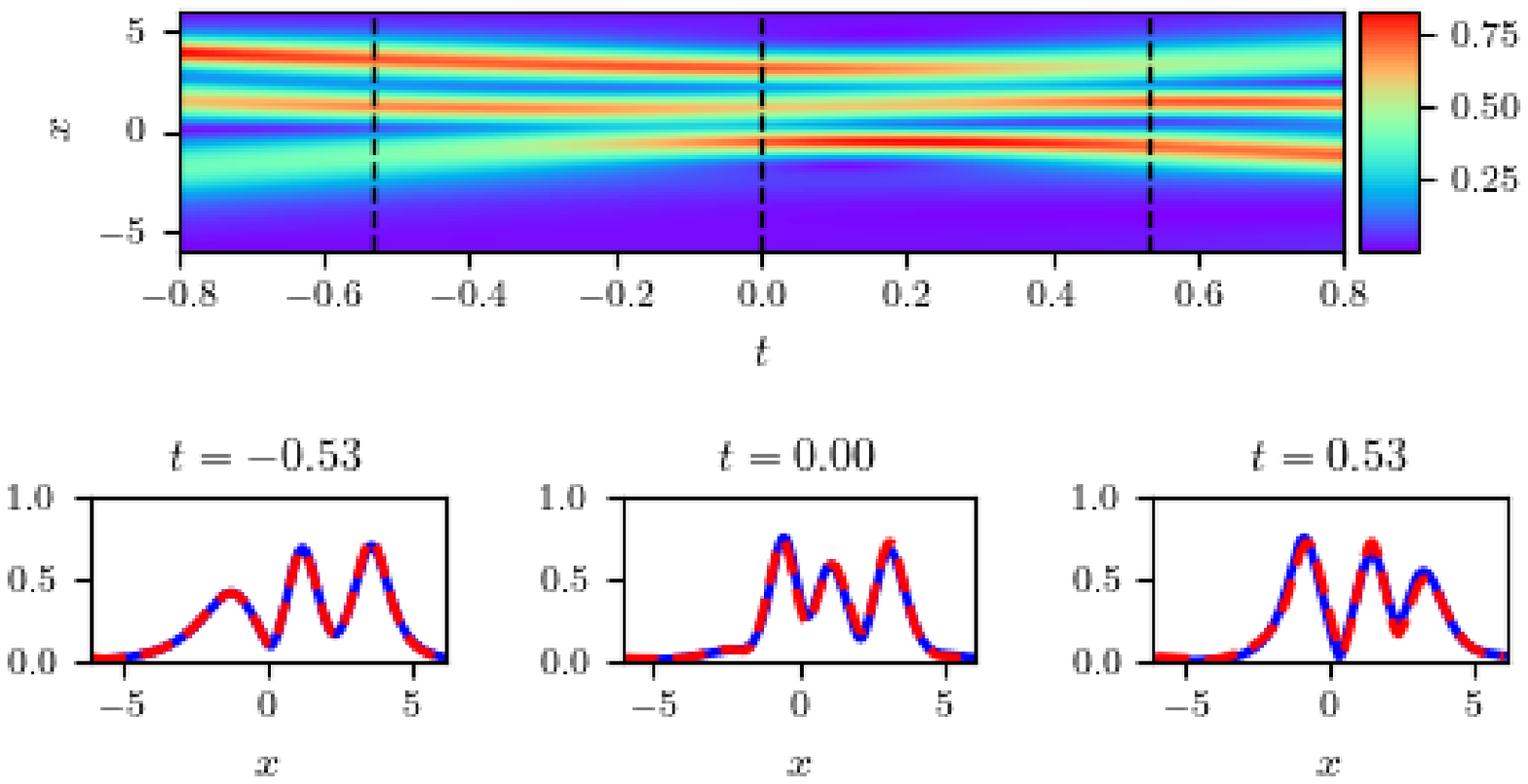}\hspace{1cm}
	\includegraphics[scale=0.4]{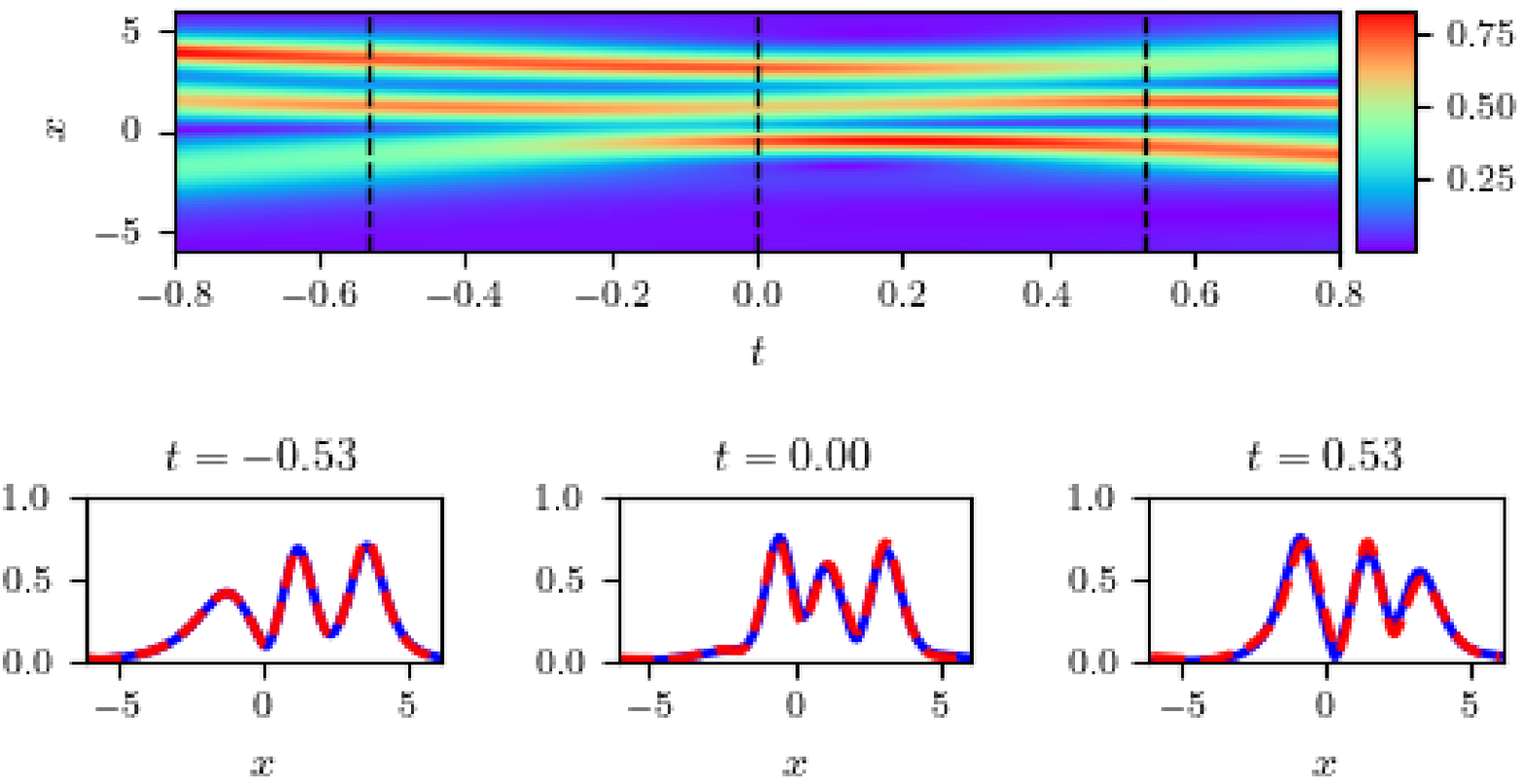}
	{\footnotesize\hspace{3cm} (c) Tractional PINN solution $ |h_{1}(x,t)| $ \hspace{3cm}(d) Tractional PINN solution $ |h_{2}(x,t)| $ } \\
	\vspace{5mm}
	\includegraphics[scale=0.4]{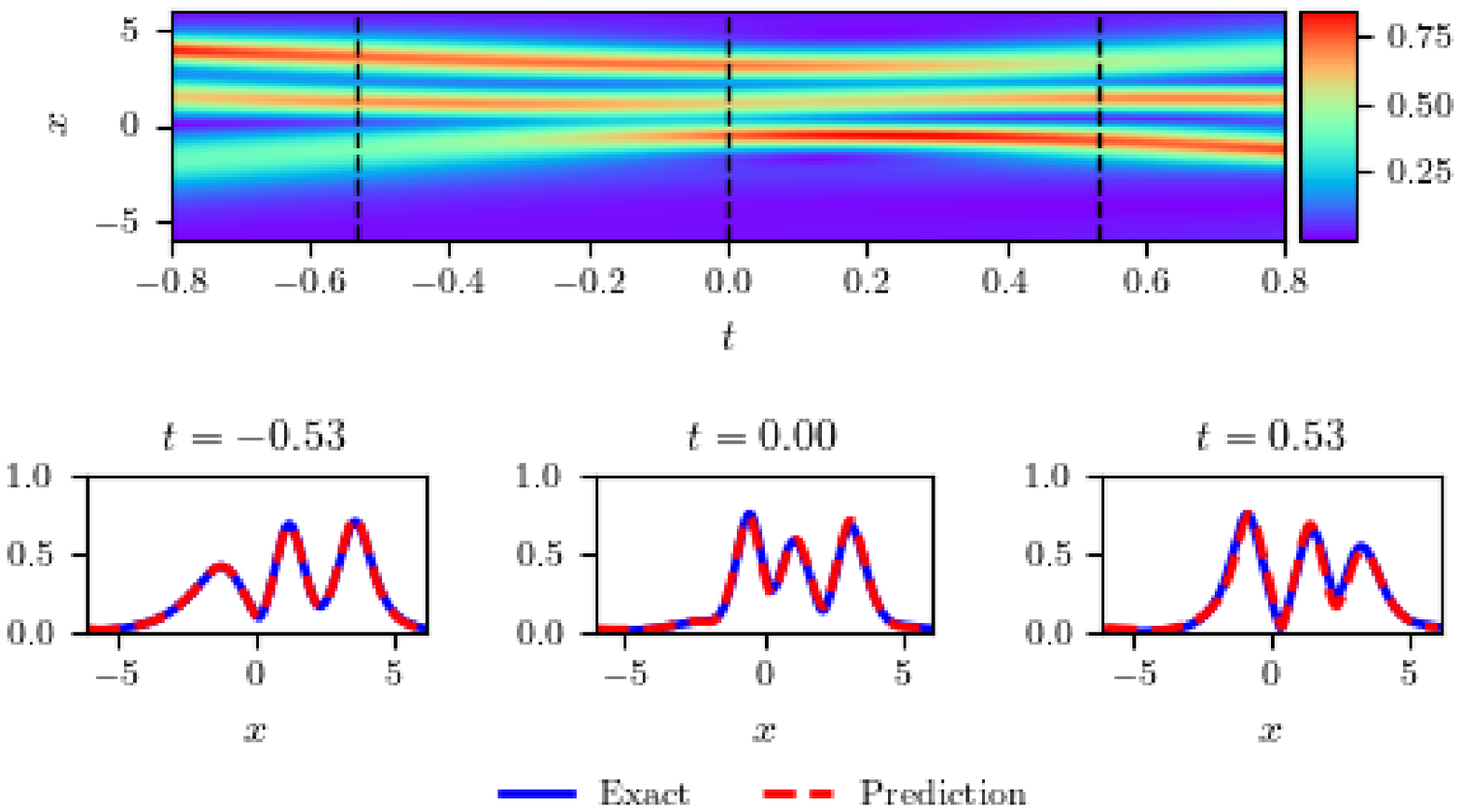}\hspace{1cm}
	\includegraphics[scale=0.4]{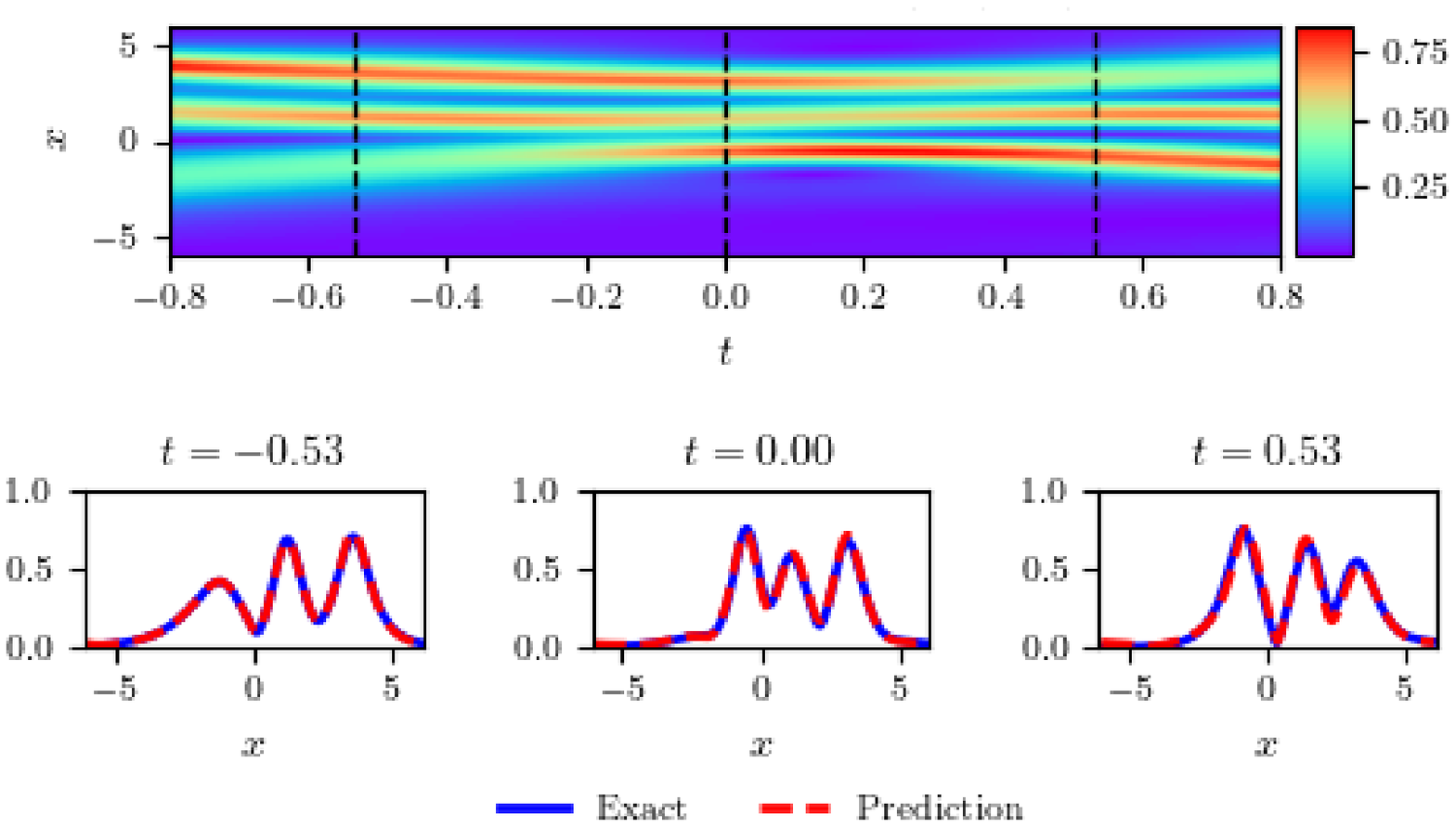}
	{\footnotesize\hspace{3.5cm} (e) RAR-PINN solution $ |h_{1}(x,t)| $ \hspace{3.5cm}(f) RAR-PINN solution $ |h_{2}(x,t)| $ }  \\
	\caption{The vector three-soliton solution $ h_{1}(x,t) $ and $ h_{2}(x,t) $: Comparison among the exact solutions (panels (a) and (b)), predictive solutions via TPINN (panels (c) and (d)) and predictive results by RAR-PINN (panels (e) and (f)) at different propagation times $ t=-0.53, 0.0 $ and 0.53. The left and right graphs represent the density diagrams and the profiles at different times for the vector three-soliton solutions $ h_{1}(x,t) $ and $ h_{2}(x,t) $, respectively.
	}
	\label{figure11}
\end{figure}

\begin{figure}[htbp]
 	\centering
 	\includegraphics[scale=0.06]{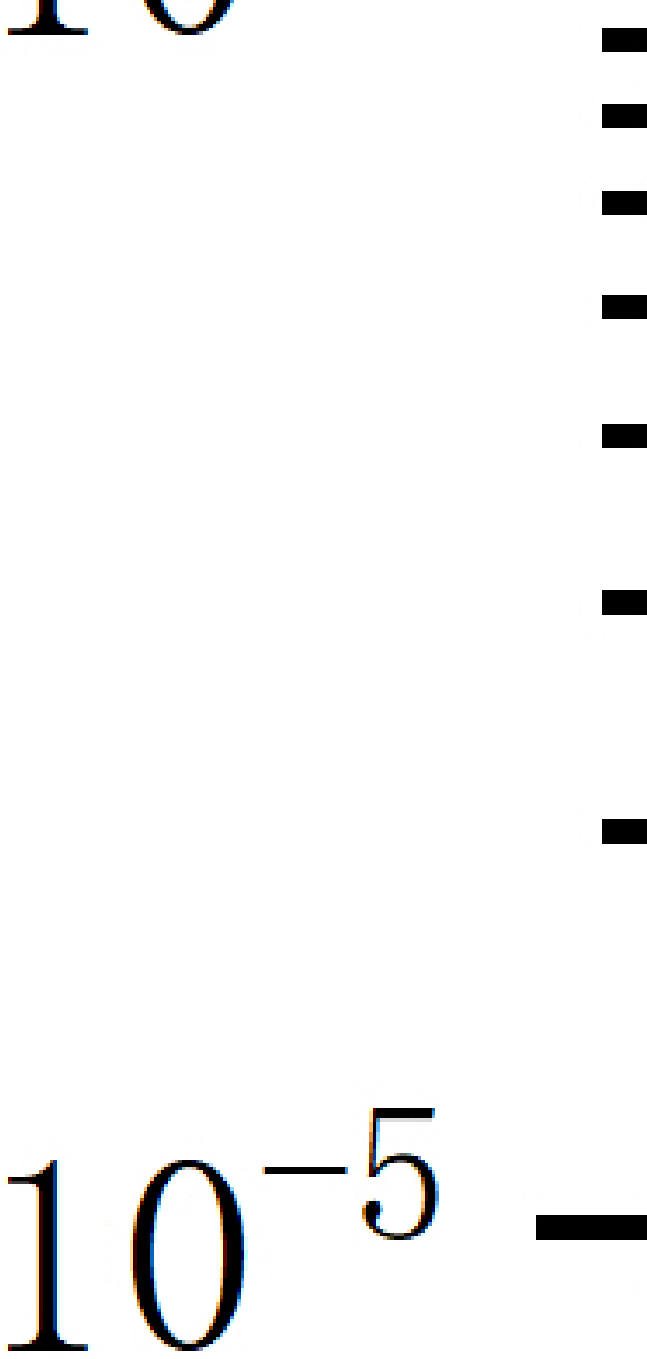}\hspace{1cm}
 	\includegraphics[scale=0.06]{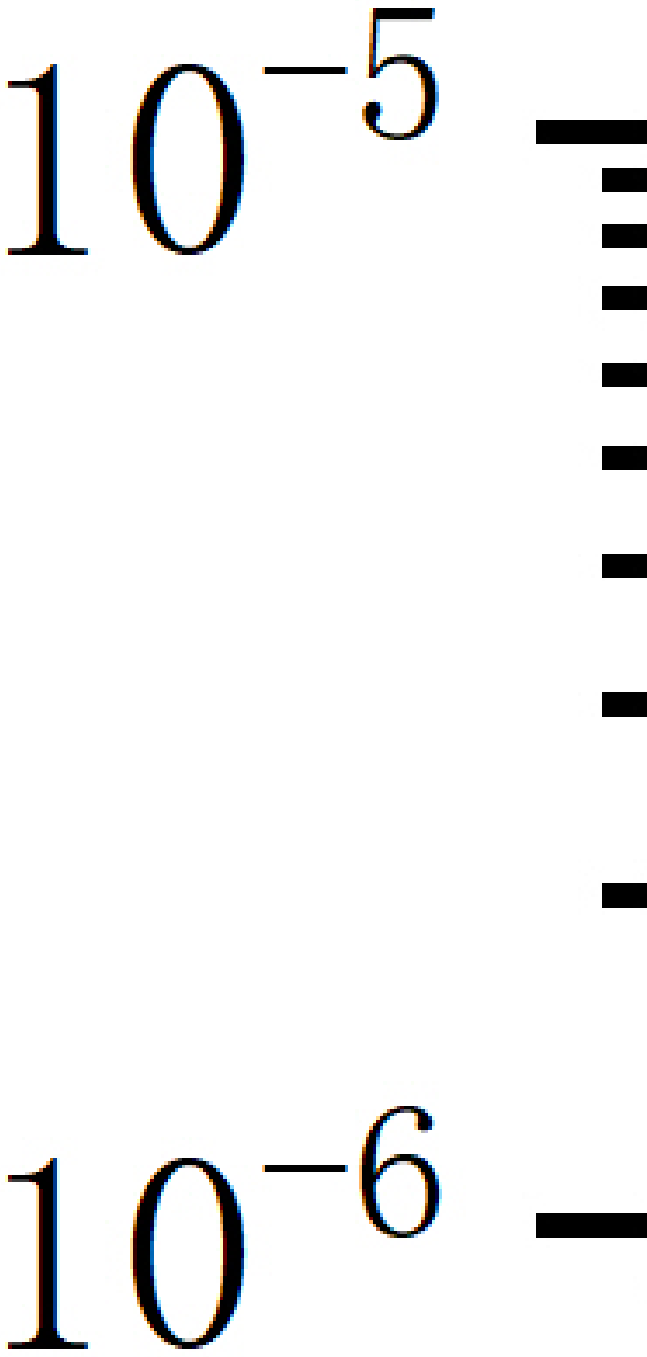}
 	{\footnotesize\hspace{9cm} (a) \hspace{8cm}(b) } \\
 	\includegraphics[scale=0.06]{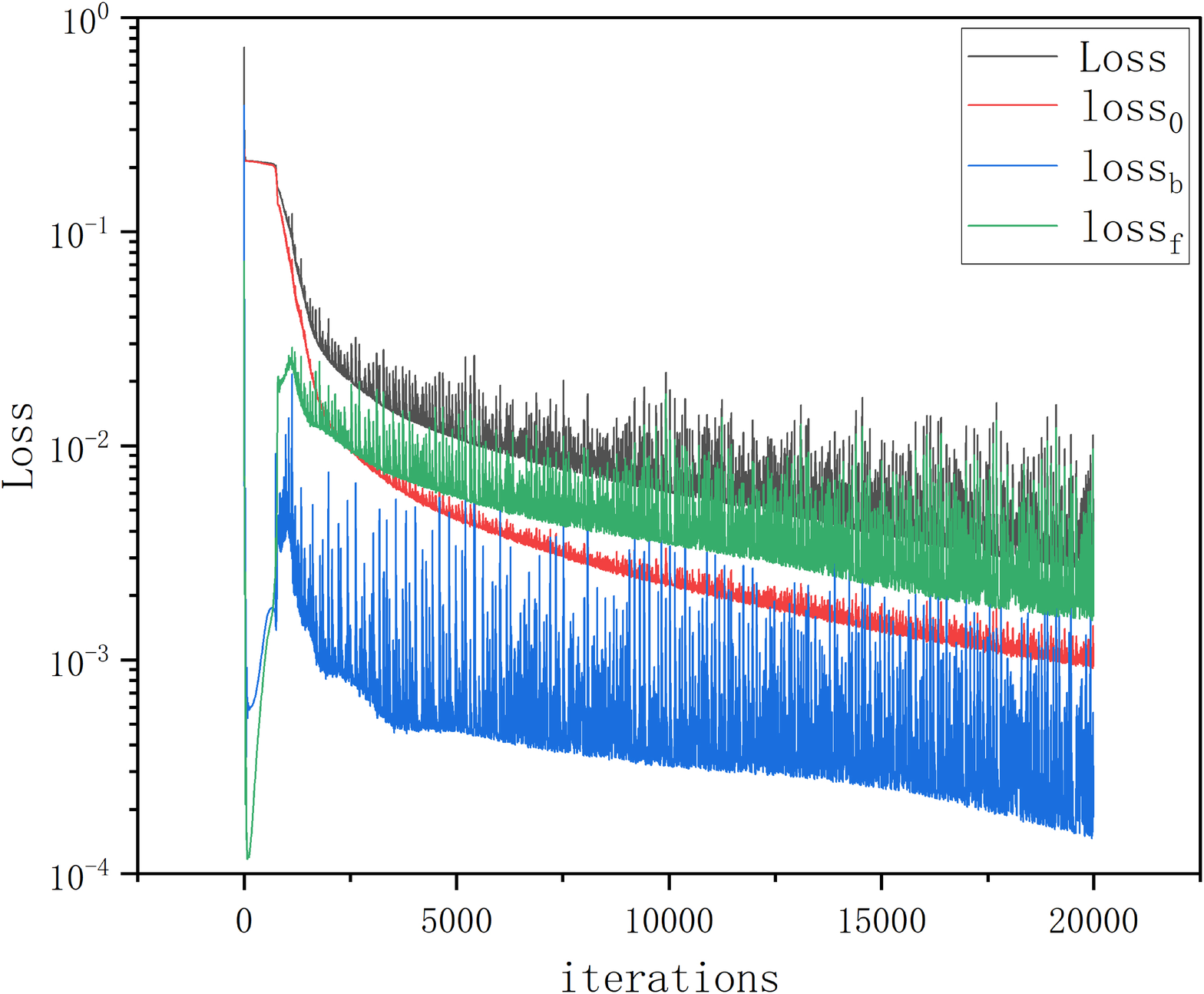}\hspace{1cm}
 	\includegraphics[scale=0.06]{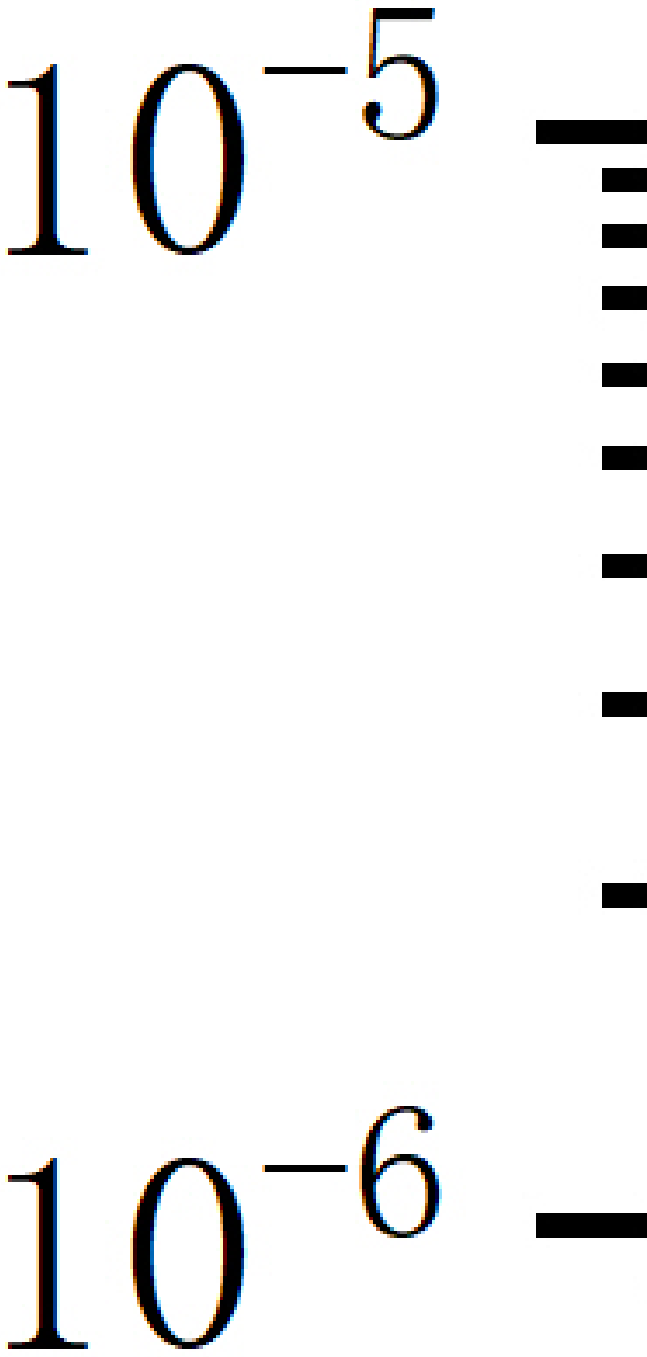}
 	{\footnotesize\hspace{9cm} (c) \hspace{8cm}(d) } \\  	
 	\caption{The loss function curve figures of vector three-soliton solutions $ h_{1}(x,t) $ and $ h_{2}(x,t) $ arising from (a) TPINN algorithm with the 20000 iterations Adam optimization; (b) TPINN algorithm with the 9893 iterations L-BFGS optimization; similarly, (c) RAR-PINN algorithm with the 20000 iterations Adam optimization; (d) RAR-PINN algorithm with the 11836 iterations L-BFGS optimization.}
 	\label{figure12}
 \end{figure}

 Figure \ref{figure11} shows the vector three-soliton solution given by the exact expression, TPINN algorithm and RAR-PINN algorithm, respectively. The relative $\mathbb{L}_{2}$ errors of the vector three-soliton solutions via the TPINN algorithm are $8.505647\times10^{-2}$ of $h_{1}(x,t)$ and $8.591699\times10^{-2}$ of $h_{2}(x,t)$, respectively; while the relative $\mathbb{L}_{2}$ errors via the RAR-PINN algorithm are $6.3181865\times10^{-2}$ of $h_{1}(x,t)$ and $6.3181865\times10^{-2}$ of $h_{2}(x,t)$. Here, we use Adam optimizer and L-BFGS optimizer to minimize the loss function. Figure \ref{figure12}(a) demonstrates the loss function descent curve of the TPINN model with the Adam optimizer for 20,000 iterations. Figure \ref{figure12}(b) indicates the loss function curve of the TPINN model with the L-BFGS optimizer for 9893 iterations. Besides, we can find the initial error $loss_{0}=1.9974363 \times10^{-5}$, the boundary error $ loss_{b}=5.9046424\times10^{-6} $, the residual $loss_{f}=5.5316887\times10^{-5}$, the total error $Loss=8.119589\times10^{-5}$. Figure \ref{figure12}(c) demonstrates the loss function curve of the RAR-PINN model via the Adam optimizer for 20,000 iterations. Figure \ref{figure12}(d) illustrates the loss function curve of the RAR-PINN model with the L-BFGS optimizer for 11836 iterations. Furthermore, we can observe the initial error $ loss_{0}=1.415507\times10^{-5}$, the boundary error $loss_{b}=3.168202\times10^{-6}$, the residual $loss_{f}=4.163423\times10^{-5}$, and the total error $Loss=5.895750\times10^{-5}$. On the whole, the RAR-PINN model is more efficient than the TPINN model in both relative and absolute errors.

\begin{figure}[htbp]
	\centering
	\includegraphics[scale=0.06]{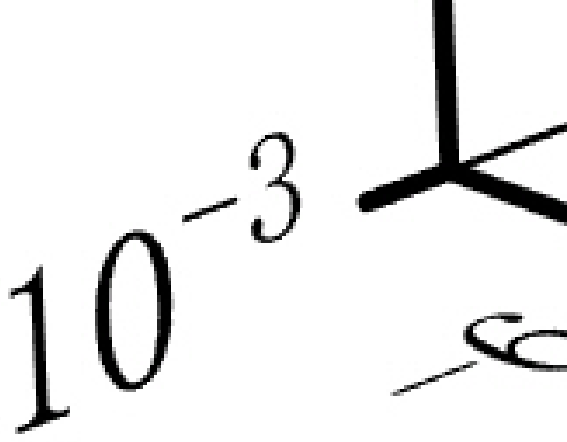}\hspace{1cm}
	\includegraphics[scale=0.06]{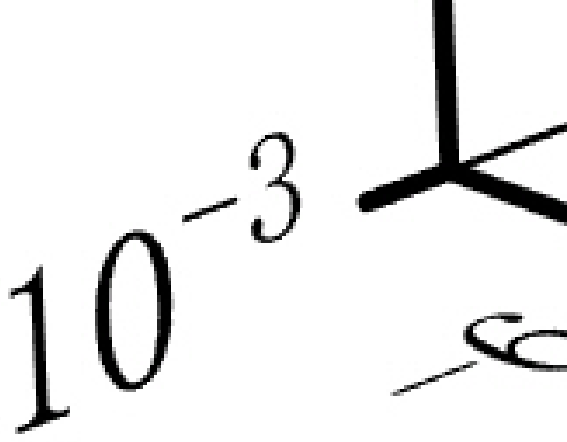}
	{\footnotesize\hspace{9cm} (a) \hspace{8cm}(b) } 	
	\caption{Three-dimensional stereogram of residual error: (a) three-dimensional stereogram residual error without adding adaptive points; (b) three-dimensional stereogram residual error after adding adaptive points.
	}
	\label{figure13}
\end{figure}

In order to further illustrate the effectiveness of RAR-PINN, three-dimensional diagram of the residual error are shown. Figure \ref{figure13} (a) and (b) show the residual errors of the three-dimensional stereogram before and after adding the adaptive points, respectively. By comparing Figure \ref{figure13} (a) and  Figure \ref{figure13} (b), we can see that the residual can reach 0.85 before adding the adaptive points and can be reduced to less than 0.45 after adding the adaptive points. Through the comparison of the above experimental results, the RAR-PINN model is more efficient for solving vector three-soliton solutions of Eq.~\eqref{equation}.

\section{Data-driven parameters discovery of the CGNLS equation}

The deep learning method can tackle the data-driven parameters discovery of the $N$-coupled nonlinear equation. For convenience, we use the CGNLS equation as an example to illustrate. In this part, there are some unknown parameters in System (1), but some extra information is known on these points besides the CGNLS equation. We will consider the problem of data-driven discovery of PDEs \cite{raissi2017machine,raissi2018hidden,rudy2017data}. The CGNLS equation with unknown parameters is as follows
\begin{equation}
	\left\{
	\begin{aligned}
		&ih_{1t}+\lambda_{1}h_{1xx}+\lambda_{2}( \alpha |h_{1}|^2+\beta|h_{2}|^2+ \gamma h_{1}h_{2}^{*}+ \gamma ^{*}h_{2}h_{1}^{*})h_{1}=0,\\
		&ih_{2t}+\lambda_{3}h_{2xx}+\lambda_{4}(\alpha |h_{1}|^2+\beta |h_{2}|^2+\gamma h_{1}h_{2}^{*}+\gamma ^{*}h_{2}h_{1}^{*})h_{2}=0,
	\end{aligned}
	\right.
\end{equation}
where $h_{k}=u_{k}+v_{k}i$ ($k=1,2$; $u_{k}$ and $v_{k}$ are the real and imaginary parts, respectively) is a complex field. Besides,  coefficients $\lambda_{1}, \lambda_{2}, \lambda_{3} $ and $ \lambda_{4}$ are unknown parameters that need to be trained. With the parameter choices $\alpha=\beta=\gamma=1$, let us consider the  physical model given explicitly by
\begin{equation}
	\left\{
	\begin{aligned}
		&f_{1}:=ih_{1t}+\lambda_{1}h_{1xx}+\lambda_{2}(|h_{1}|^2+|h_{2}|^2+ h_{1}h_{2}^{*})h_{1},\\
		&f_{2}:=ih_{2t}+\lambda_{3}h_{2xx}+\lambda_{4}(|h_{1}|^2+|h_{2}|^2+h_{1}h_{2}^{*})h_{2},\label{in}
	\end{aligned}
	\right.
\end{equation}
which contain real and imaginary parts, thus
\begin{equation}
	\left\{
	\begin{aligned}
		&f_{1u}:=v_{1t} -\lambda_{1}u_{1xx} - \lambda_{2}(u_{1}^2 + v_{1}^2 + u_{2}^2 +v_{2}^2+2(u_{1}u_{2}+v_{1}v_{2}))u_{1},\\
		&f_{1v}:=u_{1t}+\lambda_{1}v_{1xx}+\lambda_{2}(u_{1}^2 + v_{1}^2 + u_{2}^2 +v_{2}^2+2(u_{1}u_{2}+v_{1}v_{2}))v_{1},\\
		&f_{2u}:=v_{2t} -\lambda_{3}u_{2xx} - \lambda_{4}(u_{1}^2 + v_{1}^2 + u_{2}^2 +v_{2}^2+2(u_{1}u_{2}+v_{1}v_{2}))u_{2},\\
		&f_{2v}:=u_{2t}+\lambda_{3}v_{2xx}+\lambda_{4}(u_{1}^2 + v_{1}^2 + u_{2}^2 +v_{2}^2+2(u_{1}u_{2}+v_{1}v_{2}))v_{2}.
	\end{aligned}
	\right.
\end{equation}
We determine the values of the unknown parameters $ \lambda_{1}, \lambda_{2}, \lambda_{3} $ and $ \lambda_{4} $ by minimizing the MSE loss function, which is defined as follows
\begin{equation}
	MSE = MSE_{p}+MSE_{f},
\end{equation}
where
\begin{equation}
	\begin{split}
		\begin{aligned}
		MSE_{p}=&\frac{1}{N_{u}}\sum_{i=1}^{N_{u}}({|\hat{u}_{1}}(x^{i},t^{i})-u_{1}^i|^2+|\hat{v}_{1}(x^{i},t^{i})-v_{1}^i|^2\\
		&+|\hat{u}_{2}(x^{i},t^{i})-u_{2}^i|^2+|\hat{v}_{2}(x^{i},t^{i})-v_{2}^i|^2),\\
		MSE_{f}=&\frac{1}{N_{u}}\sum_{i=1}^{N_{u}}(|f_{1u}(x^i,t^i)|^2+|f_{1v}(x^i,t^i)|^2+|f_{2u}(x^i,t^i)|^2+|f_{2v}(x^i,t^i)|^2).
		\end{aligned}
	\end{split}
\end{equation}
The vector-soliton solutions are used as the data set. Using the pseudo-spectral method, the solutions of the Eq.~\eqref{in} are discretized into $[300\times201]$ data points. The exact solutions for the Eq.~\eqref{in} with parameters $\lambda_{1}=1, \lambda_{2}=2, \lambda_{3}=1$ and $\lambda_{4}=2$ are considered to make sampling data set with space-time region $(x,t)\in[-10,10]\times[-2,2]$.

To illustrate the effectiveness of our algorithm, we have created a training data-set by randomly generating $ N_{u} = 5,000 $ points across the entire domain from the exact solution corresponding to $\lambda_{1}=1, \lambda_{2}=2, \lambda_{3}=1$ and $\lambda_{4}=2$. During the training process, the RAR-PINN algorithm is used in order to add more points where the error is large. We add 10 more new points adaptively via RAR with $m=5$ and $\varepsilon_{0}=0.15$, and then the total number of data points is $N=5,010$. Table \ref{tab1} indicates that the RAR-PINN algorithm is able to correctly identify the unknown parameters with very high accuracy even when the training data are corrupted with noise. Specifically, The predictions remain robust even when the training data are corrupted with $1\%$ and $3\%$ uncorrelated Gaussian noise.

\begin{table}
	\centering
	\fontsize{10}{15}\selectfont    
	\caption{Identified obtained by learning $\lambda_{1},\lambda_{2},\lambda_{3}$ and $\lambda_{4}$}
	\begin{tabular}{ | c | c |}
		\hline		
		\multirowcell{2}{Correct CGNLS equation}
		&$ih_{1t}+h_{1xx}+2(|h_{1}|^2+|h_{2}|^2+h_{1}h_{2}^{*})h_{1}=0$\\
		&$ih_{2t}+h_{2xx}+2(|h_{1}|^2+|h_{2}|^2+h_{1}h_{2}^{*})h_{2}=0$\\
		\hline
		\multirowcell{2}{Identified CGNLS equation (clean data)}
		&$ih_{1t}+1.00012h_{1xx}+2.00270(|h_{1}|^2+|h_{2}|^2+h_{1}h_{2}^{*})h_{1}=0$\\
		&$ih_{2t}+1.00014h_{2xx}+2.00359(|h_{1}|^2+|h_{2}|^2+h_{1}h_{2}^{*})h_{2}=0$\\
		\hline
		\multirowcell{2}{Identified CGNLS equation ( $ 1 \% $ noise)}
		&$ih_{1t}+0.99912h_{1xx}+2.00500(|h_{1}|^2+|h_{2}|^2+h_{1}h_{2}^{*})h_{1}=0$\\
		&$ih_{2t}+1.00057h_{2xx}+2.00436(|h_{1}|^2+|h_{2}|^2+h_{1}h_{2}^{*})h_{2}=0$\\
		\hline
		\multirowcell{2}{Identified CGNLS equation ( $ 3 \% $ noise)}
		&$ih_{1t}+1.01020h_{1xx}+2.01574(|h_{1}|^2+|h_{2}|^2+h_{1}h_{2}^{*})h_{1}=0$\\
		&$ih_{2t}+0.99996h_{2xx}+2.01110(|h_{1}|^2+|h_{2}|^2+h_{1}h_{2}^{*})h_{2}=0$\\
		\hline
	\end{tabular}
	\label{tab1}
\end{table}

To further check the stability of our algorithm, we conduct a series of experiments with different amounts of training data and noise corruption levels. The results are summarized in Tables \ref{tab2} and \ref{tab3}. The key point here is that this method seems to be very robust to the level of noise in the data, producing reasonable recognition accuracy even for noise corruptions of up to $12\%$. We also observe from Table~\ref{tab2} that the error decreases as the total number of training points increases, while the error is non-monotonic along with the change of the noise levels as seen in Table \ref{tab3}. Different factors related to the equations themselves as well as the particular neural network setup may be responsible for this variation.

\begin{table}
	\centering
	\fontsize{10}{15}\selectfont    
	\caption{Error in the identified parameters $ \lambda_{1}, \lambda_{2}, \lambda_{3} $, and $ \lambda_{4} $ for different number of training data}
	\begin{tabular}{c|ccccc}
		\toprule
		\diagbox [width=8em,trim=r] {Parameters}{$ N $}
		&  1015  &  2010  &  3008   &  4009  &  5010  \\
		\hline
		$ \lambda_{1} $ & 1.143306 & 0.0043837 & 0.0021511 & 0.0016654 & 0.0001152 \\
		$ \lambda_{2} $ & 0.9430671 & 0.0022939 & 0.0007388 & 0.0038375 & 0.0013522 \\
		$ \lambda_{3} $ & 1.0143343 & 0.0017641 & 0.0016935 & 0.0016654 & 0.0001363 \\
		$ \lambda_{4} $ & 0.9438857 & 0.0017554 & 0.0001838 & 0.0038375 & 0.0017950 \\
		\bottomrule
	\end{tabular}\vspace{0.5cm}
	\label{tab2}
\end{table}
\begin{table}
	\centering
	\fontsize{10}{15}\selectfont    
	\caption{Error in the identified parameters $ \lambda_{1}, \lambda_{2}, \lambda_{3} $, and $ \lambda_{4} $ by different noise levels}
	\begin{tabular}{c|ccccc}
		\toprule
		\diagbox [width=8em,trim=r] {Parameters}{Noise}
		& $ 1 \% $ & $ 3 \% $ & $ 6 \% $  & $ 9 \% $ & $ 12 \% $ \\
		\hline
		$ \lambda_{1} $ & 0.0008779 & 0.0101990 & 0.0059060 & 0.0009057 & 0.0073389 \\
		$ \lambda_{2} $ & 0.0025021 & 0.0078678 & 0.0042040 & 0.0010816 & 0.0051344\\
		$ \lambda_{3} $ & 0.0005662 & 0.0000419 & 0.0024133 & 0.0015854 & 0.0059356 \\
		$ \lambda_{4} $ & 0.0021807 & 0.0055518 & 0.0020984 & 0.0016636 & 0.0051588 \\
		\bottomrule
	\end{tabular}\vspace{0.5cm}
	\label{tab3}
\end{table}

\section{Discussion and conclusion}

In conclusion, we have investigated the CGNLS equation via the RAR-PINN method. By improving the distribution of residual points during the training process, the RAR-PINN method increases the training efficiency. We have used RAR-PINN and TPINN methods to predict the vector one-, two-, and three-soliton solutions of the CGNLS equation, respectively, and then compared the experiment results. It has been revealed that the RAR-PINN method can present higher accuracy for the same number of training points. In particular, the RAR-PINN method still achieves good results for shape-changing soliton solutions of the CGNLS equation. After a series of experiments, we can conclude that when solving the coupled nonlinear equations, the RAR-PINN method can well simulate vector-soliton interactions for both elastic and inelastic collisions. In addition, the data-driven parameters discovery of the CGNLS equation has been performed. We have added the noise to judge the stability of the neural networks. The prediction error increases as the number of sampling points decreases. These changes of error are within the controllable range, which is enough to prove the excellent performance of the RAR-PINN method. Such method can be applied into other coupled nonlinear systems to discover the dynamic properties of nonlinear waves.
\section*{Acknowledgments}
This work was partially supported by the Natural Science Foundation of Beijing Municipality (Grant No. 1212007), by the National
Natural Science Foundation of China (Grant No. 11705284), by the Fundamental Research Funds of the Central Universities (Grant
No. 2020MS043), and by the Science Foundations of China University of Petroleum, Beijing (Grant Nos. 2462020YXZZ004 and 2462020XKJS02).
\bibliography{reference}

\end{document}